\documentclass[a4paper,11pt,oneside]{article}
\usepackage{lineno}
\usepackage{hyperref}

\usepackage{amsmath,amsthm,amsfonts,amssymb,amscd}
\usepackage{graphicx}
\usepackage{mathtools}
\usepackage{hyperref}
\usepackage{pgf,tikz}
\usepackage{mathrsfs}
\usepackage{float}
\usepackage{csvsimple}
\usepackage{enumerate}
\usepackage{booktabs}
\usepackage{pgfplots}
\pgfplotsset{compat=1.15}
\usepackage{pgfplotstable}
\usepackage[noend]{algpseudocode}
\usepackage{graphicx,bm,xcolor}
\usepackage{algorithm}
\usepackage{pdfpages}
\usepackage[noend]{algpseudocode}
\usepackage{stmaryrd}
\usepackage{hyperref}
\usetikzlibrary{arrows}
\usepackage{caption}
\usepackage{t1enc}
\usepackage[polish,german,english]{babel}
\usepackage{verbatim} 
\usepackage{xcolor}

\usepackage{url}

\newif\ifMAKEPICS
\MAKEPICSfalse

\ifMAKEPICS
\usepackage[cleanup,subfolder]{gnuplottex}
\usepackage{xparse}

\ExplSyntaxOn
\DeclareExpandableDocumentCommand{\convertlen}{ O{cm} m }
{
	\dim_to_decimal_in_unit:nn { #2 } { 1 #1 } cm
}
\ExplSyntaxOff
\fi

\newcommand\plot[1]{\let\frame\relax
	\frame{\includegraphics[clip,trim=0 220 0 60,width=8cm]{#1}}}

\usepackage{authblk}

\theoremstyle{plain} 
\newtheorem{theorem}{Theorem}[section]

\newtheorem{Corollary}[theorem]{Corollary}
\newtheorem{Remark}[theorem]{Remark}
\newtheorem{Definition}[theorem]{Definition}

\newtheorem{Problem}[theorem]{Problem}

\theoremstyle{definition} %
\newtheorem{assumption}{Assumption}

\theoremstyle{remark} %

\usepackage{newunicodechar}

\usepackage{xcolor}
\newcommand{\bernhard}[1]{{{\textcolor{black}{#1}}}}
\newcommand{\twick}[1]{{\textcolor{black}{#1}}}
\newcommand{\Sven}[1]{{\textcolor{black}{#1}}}

\begin{document}

		\title{
			{{Mathematical} modeling and {numerical} multigoal-oriented a posteriori error control and adaptivity for a stationary, nonlinear, coupled
				flow temperature model with temperature dependent density}
		}
		
		\author[1,3]{S. Beuchler}
		\author[2,3]{A. Demircan}
		\author[1,3]{B. Endtmayer}
		\author[2,3]{U. Morgner}
		\author[1,3]{T. Wick}
		
\affil[1]{Cluster of Excellence PhoenixD (Photonics, Optics, and
	Engineering -- Innovation Across Disciplines), Leibniz Universit\"at Hannover, Germany}
\affil[2]{Leibniz Universit\"at Hannover, Institute of Quantum Optics, Welfengarten 1, 30167 Hannover, Germany}
\affil[3]{ Leibniz Universit\"at Hannover, Institute of Applied Mathematics , Welfengarten 1, 30167 Hannover, Germany}

\date{}

\maketitle
		
		\begin{abstract}
			
			In this work, we develop adaptive schemes using goal-oriented 
			error control for a highly nonlinear flow temperature model with temperature dependent density.
			The dual-weighted residual method for computing error indicators to steer mesh 
			refinement and solver control is employed. The error indicators are used to 
			employ adaptive algorithms, which are substantiated with several numerical tests. 
			Therein, error reductions and effectivity indices are consulted to establish 
			the robustness and efficiency of our framework.
			
		\end{abstract}


	\section{Introduction}
	This work is devoted to the development of a mathematical model of the Navier-Stokes equations coupled 
	to a heat equation with temperature dependent viscosity and density. 
	The time dependent Navier-Stokes equations with density were investigated in 
	\cite{fernandez2012motivation}, while stationary Navier-Stokes equations with variable density (without temperature) were considered in \cite{mallea2019bilinear}.
	Well-known are so-called Boussinesq approximations \cite{boussinesq_review,LoBo96}.
	The resulting problem is a highly nonlinear coupled PDE (partial differential equations) system. Applications of temperature dependent viscosities 
	and densities range from ionic liquids \cite{B513231B} over laser material 
	processing \cite{OTTO201035} up to wave guide modeling \cite{Pae17,Chenetal18}. 
	{
		Related theoretical work is \cite{FeiMa06} in which long-time and large-data existence of a weak solutions were obtained with the viscosity and heat-conductivity coefficients varying with the temperature.

		Moreover, nonstationary Navier-Stokes with variable densities and variable viscosities 
		were considered in \cite{AxHeNe15}.
		Finally, \cite{Pech2020} proposed a computational algorithm for the 
		evolutionary Navier-Stokes-Fourier system with temperature dependent material properties such as 
		density, viscosity and thermal conductivity.
		To the best of our knowledge, exactly our proposed model and numerical realization has not yet been utilized in the literature.
	}
	
	Often in numerical simulations,
	not the entire numerical solution is of interest, but certain quantities 
	of interest (QoI). Accurate evaluations are obtained with the help 
	of goal-oriented error control with the dual-weighted residual method (DWR)~\cite{BeRa01} (see also \cite{AinsworthOden:2000,BaRa03,Od18}), which serve as basis for 
	mesh adaptivity and balancing discretization and nonlinear iteration errors. 
	The latter ensures that solver tolerances are automatically adjusted according to the current discretization accuracy of the QoIs. As we have 
	a multiphysics situation at hand, multiple QoIs might be of interest simultaneously. This is addressed by multigoal-oriented error control using again 
	the DWR method. Early work was carried out in \cite{Ha08}. Significant 
	improvements from the algorithmic and theoretical level were done by us 

	in prior work, 
        namely a partition-of-unity localization (based on \cite{RiWi15_dwr})
        for multigoal-oriented error control in \cite{EndtWi17}, balancing 
        discretization and nonlinear iteration errors for multigoal in \cite{EnLaWi18},
        efficiency and reliability estimates for goal-oriented error control 
        using a saturation assumption in \cite{EndtLaWi20}, and 
        smart algorithms switching between solving high-order 
        problems or using interpolations in \cite{EndtLaWi21_smart}.
{        Additionally, efficiency and reliability estimates for arbitrary interpolations are established in \cite{EndtLaWi21_smart}.}

	Other conceptional improvements in multigoal
	error estimation were done in \cite{BruZhuZwie16,PARDO20101953,BREVIS2021186}.
	Recent applications in stationary multiphysics 
	problems were undertaken in fluid-structure interaction \cite{AhEndtSteiWi22} and 
	the Boussinesq equations \cite{BeuEndtLaWi23}. 
	
	The objective in this work is a nontrivial extension of \cite{BeuEndtLaWi23}
	by making the density variable. This results into both 
	temperature-dependent densities and viscosities. Their respective 
	relationships are described with exponential laws. The resulting 
	coupled problem is highly nonlinear. The numerical solution algorithms 
	are based on monolithic formulations in which the entire system is solved 
	all-at-once. Therein, a Newton solver is utilized for the nonlinear solution in which the Newton tolerances are chosen according to the current accuracy of the quantities of interest, i.e., balancing discretization and nonlinear iteration errors.
	This is achieved by a multigoal-oriented a posteriori error estimates 
	with adjoint problems, using the DWR method, and was for multigoal error estimation first formulated in \cite{EnLaWi18}. The adjoint needs to be explicitly derived for our coupled problem. In order to utilize the error estimators 
	for local mesh adaptivity, the error needs to be localized, which is 
	here done with a partition-of-unity (PU), proposed for DWR in \cite{RiWi15_dwr}.
	These developments result into our final adaptive algorithms.
	Our numerical examples are inspired from applications in laser material processing.
	
	The outline of this paper is as follows: In Section~\ref{sec_model_problem} 
	our model problem and discretization with the Galerkin finite element method 
	are introduced. Next, in Section~\ref{sec_DWR} our goal-oriented framework 
	is introduced. Section~\ref{Sec: Num. Example} is the main part of this paper. 
        Several numerical examples are conducted and investigated with respect to error reductions and effectivity indices.
	In most of the experiments \twick{point-like sources} of the heat equation is chosen. The coarse mesh of the approximation requires some a-priori refinement which is explained in Subsection~\ref{Pre-refinement for laser}.
	This sections contains also comparisons of the simpler model investigated in \cite{BeuEndtLaWi23} with the new one. 
	Also computations on a non-polygonal geometry are presented.
	Finally, our work is 
	summarized in Section~\ref{Sec; Concl}.

	\section{The model problem and discretization}
	\label{sec_model_problem}
	\subsection{The model problem}
	Let $\Omega \subset \mathbb{R}^d$, with $d \in \{2, 3\}$, be a bounded Lipschitz domain.
	As model problem, we consider the Navier-Stokes equation with temperature dependent viscosity and density.
	The corresponding problem looks as follows: Find the velocity field $v$,the pressure $p$ and the temperature $\theta$ such that
	\begin{equation}
	\begin{aligned}
	-div\left(\nu(\theta)\rho(\theta)(\nabla v + \nabla v^T\right) + (\rho(\theta) v \cdot \nabla) v  - \nabla p - g \rho(\theta)&=0,  \nonumber \\
	-div(\rho(\theta)v)&=0, \nonumber\\
	-div(k \nabla \theta) + v\cdot \nabla \theta & = f_{E},\label{eq:stationary_navierstokes}\\
	\rho(\theta) &= \rho_0e^{-\int_{\theta_0}^{\theta}\alpha(\hat{\theta})\text{ d}\hat{\theta}},\nonumber\\
	\nu(\theta) &= \nu_0 e^{\frac{E_A}{R \theta}},\nonumber
	\end{aligned} 
	\end{equation}
	\text{in $\Omega$.} In this paper, we use Dirichlet boundary conditions for velocities and temperature only.
	More precisely
	\begin{equation*}
	\begin{aligned}
	v&=0  \qquad \text{ on } \partial\Omega, \\
	\theta&=\theta_0 \qquad  \text{ on }\partial \Omega,\\
	\end{aligned}
	\end{equation*}
	where $k>0$ is the thermal conductivity and $\theta_0>0$. 
	Of course the conditions  $v=0 \text{ on }\partial \Omega$ and $\theta=\theta_0 \text{ on }\partial \Omega$ have to be understood as traces.
	However, also other boundary conditions do not restrict the result of this work.
	The viscosity $\nu (\theta)$ is modeled by the Arrhenius equation; see \cite{Arrhenius+1889+226+248,Arrhenius+1889+96+116} for chemical reactions and \cite{de1913relation,raman1923theory,andrade1934lviii,ward1937viscosity,haj2014contribution} for viscosity. Here, the parameters $E_A$ and $\nu_0$ are material parameters and $R=8.31446261815324$ is the universal gas constant.
	The equation for the density is derived from the volume equations in \cite{DUBROVINSKY20021}, i.e. $$V(\theta) = V_0e^{\int_{\theta_0}^{\theta}\alpha(\hat{\theta})\text{ d}\hat{\theta}}.$$	
	Since the density depends on the temperature, the part $g \rho(\theta)$ is not necessarily a gradient field and therefore might change the flow. For constant density, this term just influences the pressure.
	Furthermore let 
	\begin{align*}
	V_{v,0}&:= [H_0^1(\Omega)]^d , \quad
	V_p:=L^2_0(\Omega), \quad\textrm{and}\quad
	V_{\theta,0}:=H_0^1(\Omega).
	\end{align*}
	For our model problem, we use the following weak form: 
	\begin{Problem}
		Find $U:=(v,p,\theta) \in \{0,0,\theta_0 \} 
		{+ \{V_{v,0}\times V_p \times V_{\theta,0}\}} =:U^D+V$ such that
		\begin{equation} \label{weak form}
		\begin{aligned}
		&(\nu(\theta)\rho(\theta)(\nabla v + \nabla v^T), \nabla \psi_v ) + (\rho(\theta)\nabla v \cdot v, \psi_v)+(p, div(\psi_v))-(\rho(\theta)g,\psi_v)
		\\+&(div(v)\rho(\theta)-\rho'(\theta)\nabla \theta \cdot v, \psi_p)
		\\+&(k\nabla\theta,\nabla \psi_\theta)+ (v \cdot \nabla \theta, \psi_\theta)-(f_E,\psi_\theta)=0
		\end{aligned}
		\end{equation}
		holds for all $\Psi:=(\psi_v,\psi_p,\psi_\theta) \in V$. Additionally, we note that the temperature $\theta$ has to be positive.
		As usual, we reformulate the weak form \eqref{weak form} as the operator equation 
		$$A(U)=0 \qquad \text { in } V^*,$$
		where $V^*$ denotes the dual space of $V$.
	\end{Problem}

	\subsection{Discretization and numerical solution}
	For the discretization, we decompose our domain into hypercubal elements $K \in \mathcal{T}_h$, where 
	$\bigcup_{K \in \mathcal{T}_h}\overline{K}=\overline{\Omega}$ and for all $K,K' \in \mathcal{T}_h:$  $K \cap K' = \emptyset$ if and only if $K\not=K'$. 
	Let $\mathbb{Q}_k(K)$ be the tensor product space of polynomials of degree $k$ on the reference domain $\hat{K}$ and  $$Q_k:=\{v_h \in \mathcal{C}(\Omega): \psi_h \in \mathbb{Q}_k(K)\quad \forall K \in \mathcal{T}_h \},$$
	the finite element space of continuous finite elements of degree $k$ on  $\mathcal{T}_h$.
	This leads the usual standard continuous tensor product finite elements on hypercubes.
	Finally we consider our finite element space for the model problem as $V_h:= Q_k^d \times Q_m \times Q_n$.		
	Finally the discretized model problem is given by:
	Find $U_h:=(v_h,p_h,\theta_h) \in U^D+V_h$ such that
	\begin{equation} \label{discretized weak form}
	A(U_h)(\Psi_h)=0,
	\end{equation}
	for all $\Psi_h \in V_h$.
	In the numerical experiments, we use $Q_2$ 
	finite elements for the velocity and $Q_1$ finite elements for the pressure and temperature. \bernhard{Later, we require an enriched space as in \cite{EndtWi17}.} For the enriched spaces, which will be introduced in the following section, we use $Q_4$ for the velocity and $Q_2$ for pressure and temperature.
	Once a basis is chosen this leads to a system of nonlinear equations. Here, we use Newton's method to solve these arising systems.

	\section{Goal oriented error estimation with the dual-weighted residual method}
	\label{sec_DWR}
	Let $U\in U^D+V$ be the solution of our model problem \eqref{weak form}. The solution $U$ might not be of primary interest, but in one functional evaluation $J:U^D+V \mapsto\mathbb{R}$ or even several functional evaluations $J_1, \ldots, J_N$ evaluated at $U$.
	Of course, if we know $U$, then we can compute these values. However, in general, we get the discrete solution $U_h \in U^D+ V_h$ instead of $U$, which introduces an error in the functional or functionals. To be precise, in many cases we do not even acquire $U_h$ since our model problem is non-linear, but some approximation $\tilde{U}$ of $U_h$. To estimate the error, i.e $J(U)-J(U_h)$, we use the dual weighted residual (DWR) method \cite{BeRa01,BeRa96,RanVi2013}.
	
	\subsection{Single goal}
	In this subsection, we derive the error estimation for one goal functional $J$. 
	\paragraph{The classical approach}$~$\\
	First of all, we introduce the adjoint problem: Find $Z:=(z_v,z_p,z_\theta)\in V$ such that
	\begin{equation}\label{Goee: cont. adjoint}
	A'(U)(\Psi,Z)=J'(U)(\Psi), \qquad \forall \Psi \in V,
	\end{equation}
	and for $U$ as a solution to the model problem \eqref{weak form}. Here, $A'(U)$ and $J'(U)$ describe the Fréchet-derivatives of $A$ and $J$, respectively. 
	The adjoint problem is linear, even if the model problem is highly nonlinear. However, finding $Z$ still requires us to solve a partial differential equation. Discretizing in the same way as for our model problem leads to the discretized adjoint problem: Find $Z_h \in V_h$ such that
	\begin{equation*}\label{Goee: disc. adjoint}
	A'(U_h)(\Psi_h,Z_h)=J'(U_h)(\Psi_h), \qquad \forall \Psi_h \in V_h
	\end{equation*}  
	and for $U_h$ as a solution to the discretized model problem \eqref{discretized weak form}.
	\begin{theorem}[Error representation(see \cite{RanVi2013,EnLaWi18,EnLaNeWiWo20})]\label{thm: error-identity}
		Let $\tilde{Z}\in V_h$ be an approximation to $Z_h$ and $\tilde{U} \in U^D+V_h$ be an approximation to $U_h$. Furthermore let $U$ be the solution to the model problem \eqref{weak form} and $Z$ be the solution to the adjoint problem \eqref{Goee: cont. adjoint}.
		Additionally let $A\in \mathcal{C}^3$ and $J \in \mathcal{C}^3$. Then, we have the following error identity
		\begin{align}
		\label{Goee: cont. error identity.}
		J(U)-J(\tilde{U})=\frac{1}{2}\left(\eta_p(\tilde{U})(Z-\tilde{Z})+\eta_a(\tilde{U},\tilde{Z})(U-\tilde{U})\right)+\eta_k(\tilde{U})(\tilde{Z})+\mathcal{R}^{(3)},
		\end{align}
		with
		\begin{align}
		\eta_p(\tilde{U})(\cdot):=&-A(\tilde{U})(\cdot), \label{Goee: primal estimator}\\
		\eta_a(\tilde{U},\tilde{Z}):=&J'(\tilde{U})(\cdot)-A'(\tilde{U})(\cdot, \tilde{Z}),\label{Goee: adjoint estimator}\\
		\eta_k(\tilde{U})(\tilde{Z}):=&-A(\tilde{U})(\tilde{Z}), \label{Goee: iteration error estimator}
		\end{align}
		and $\mathcal{R}^{(3)}$ as higher order term. For more information on $\mathcal{R}^{(3)}$, we refer to \cite{EndtLaWi20,Endt21}.
		\begin{proof}
			For the proof we refer to \cite{RanVi2013,EnLaWi18,EnLaNeWiWo20}.
		\end{proof}
	\end{theorem}
	Usually, the values $U-\tilde{U}$ and $Z-\tilde{Z}$ are not available. It is possible to approximate these quantities by using interpolation techniques \cite{BeRa01,BaRa03,RanVi2013,EndtLaWi21_smart}, or hierarchical approaches (by higher order finite elements or a finer mesh) as done in \cite{BeRa01,BaRa03,RiWi15_dwr,bause2021flexible,bruchhauser2020dual,VANDERZEE20112738}. \bernhard{A comparison between higher order finite elements and a finer mesh for the Navier-Stokes equation is provided in \cite{EnLaThieWienumath}.}
	
	\paragraph{An approach on discrete spaces}$~$\\
	Here, we introduce an enriched space $V_h^{(2)}$ with $V_h\subset V_h^{(2)}\subset V$ as it is done in the hierarchical approaches above. We introduce the enriched model problem: Find $U_h^{(2)} \in U^D+V_h^{(2)}$ such that
	\begin{equation}
	\label{Goee: enriched model p}
	A(U_h^{(2)})=0 \qquad \text{ in } V_h^{(2),*},
	\end{equation}
	where $V_h^{(2),*}$ is the dual space of $V_h^{(2)}$.
	Moreover, the enriched adjoint problem is given by: Find $Z_h^{(2)}\in V_h^{(2)}$ such that
	\begin{equation}\label{Goee: enriched disc. adjoint}
	A'(U_h^{(2)})(\Psi_h^{(2)},Z_h^{(2)})=J'(U_h^{(2)})(\Psi_h^{(2)}), \qquad \forall \Psi_h^{(2)} \in V_h^{(2)},
	\end{equation}
	with $U_h^{(2)}$ as the solution of the enriched model problem \eqref{Goee: enriched model p}.
	
	\begin{assumption}[Saturation assumption; see \cite{EndtLaWi20,EnLaSchaf2023}] \label{Ass: saturation}
		Let $U_h^{(2)}$ be the solution of \eqref{Goee: enriched model p}.
		Furthermore, let $\tilde{U}$ be an approximation to \eqref{weak form}.
		Then there is a constant $b_0 \in (0,1)$ and a constant $b_h \in (0,b_0)$ such that 
		\begin{align*}
		|J(U)-J(U_h^{(2)})|<b_h|J(U)-J(\tilde{U})|.
		\end{align*}
	\end{assumption}
	\begin{Corollary}
		\label{Goee: result_discrete}
		Let $A:U^D+V\mapsto V^*$ and $J:\mathcal{D}(J)\mapsto \mathbb{R}$. Furthermore, let us assume that $U \in \mathcal{D}(J)$, where $U$ solves the model problem \eqref{weak form}. Moreover, we assume that the operator $A$ and the functional $J$ are at least three times Fréchet differentiable on the discrete space $U^D+V_h^{(2)}$. In addition, let $U_h^{(2)} \in U^D+V_h^{(2)}$ be the solution of the enriched model problem \eqref{Goee: enriched model p} and $Z_h^{(2)}\in V_h^{(2)}$ be the enriched adjoint solution \eqref{Goee: enriched disc. adjoint}.
		Then for arbitrary $\tilde{U}\in U^D+V_h^{(2)},\tilde{Z}\in V_h^{(2)}$ it holds,
		\begin{align*}
		J(U)-J(\tilde{U})=J(U)-J(U_h^{(2)})+\frac{1}{2}\left(\eta_p(\tilde{U})(Z_h^{(2)}-\tilde{Z})+\eta_a(\tilde{U},\tilde{Z})(U_h^{(2)}-\tilde{U})\right)\\
		+&\eta_k(\tilde{U})(\tilde{Z})+\mathcal{R}_h^{(3),(2)},
		\end{align*}
		where $\mathcal{R}_h^{(3),(2)}$ is a higher order term in $\left(\Vert Z_h^{(2)}-\tilde{Z}\Vert + \Vert U_h^{(2)}-\tilde{U}\right) \Vert $ and $\eta_a, \eta_p$ and $\eta_k$ are defined by \eqref{Goee: primal estimator},\eqref{Goee: adjoint estimator} and \eqref{Goee: iteration error estimator} as in Theorem~\ref{thm: error-identity}.
		\begin{proof}
			The assertion follows immediately from Corollary 4.4 in \cite{EnLaSchaf2023}.
		\end{proof}
	\end{Corollary}

	\begin{Remark}
		By replacing $U$ by $U_h^{(2)}$ in $\eta_a, \eta_p$ and $\eta_k$ in Theorem~\ref{thm: error-identity}, we observe that they coincide with $\eta_a, \eta_p$ and $\eta_k$ in Corollary~\ref{Goee: result_discrete}.	
		However, in contrast to Theorem~\ref{thm: error-identity}, it is not required that $A$ and $J$ are three times Fréchet differentiable on the continuous level. In addition, the existence of the adjoint solution is not required anymore. However, the existence of the enriched adjoint solution is required instead. Additionally we need the existence of the functional value $J(U)$.
	\end{Remark}
	\begin{Definition}[Efficient and reliable error estimator for $J$]
		The error estimator $\eta$ is {efficient}, if there exists a constant $\underline{c} \in \mathbb{R}$ with $\underline{c}>0$ such that
		\begin{equation*} \label{eq: efficient}
		\underline{c}|\eta| \leq |J(U)-J(\tilde{U})|.
		\end{equation*}
		We say an error estimator $\eta$ is {reliable}, if there exists a constant $\overline{c} \in \mathbb{R}$ with $\overline{c}>0$ such that
		\begin{equation*} \label{eq: reliable}
		\overline{c}|\eta| \geq |J(U)-J(\tilde{U})|.
		\end{equation*}
	\end{Definition}
\begin{Definition}
	Following the result of Corollary \ref{Goee: result_discrete}, we define our error estimator as 
	\begin{align*}
	\eta:=\overbrace{\frac{1}{2}\left(\eta_p(\tilde{U})(Z_h^{(2)}-\tilde{Z})+\eta_a(\tilde{U},\tilde{Z})(U_h^{(2)}-\tilde{U})\right)}^{:=\eta_h}+\eta_k(\tilde{U})(\tilde{Z})+\mathcal{R}_h^{(3),(2)}.
	\end{align*}
\end{Definition}
In the work \cite{EndtLaWi20} was shown, that the part $\mathcal{R}_h^{(3),(2)}$ indeed is of higher order and therefore can be neglected.
	The part $\eta_k(\tilde{U})(\tilde{Z})$ can be used to stop the nonlinear solver as suggested in \cite{RanVi2013,EndtLaWi20}.
	The following theorem shows efficiency and reliability for the resulting error estimator.
	\begin{Corollary}[Efficiency and reliability for $\eta$]\label{Corr: efficiency and reliablity eta}
		Let us assume that Assumption~\ref{Ass: saturation} is fulfilled. Furthermore let all assumption of Corollary~\ref{Goee: result_discrete} be fulfilled. Then
		the error estimator $\eta$ is efficient and reliable with the constants $\overline{c}=\frac{1}{1-b_h}$ and $\underline{c}=\frac{1}{1+b_h}$.
		\begin{proof}
			The proof follows the same ideas as the proof in \cite{EndtLaWi20}.
			From Corollary~\ref{Goee: result_discrete} we know that
			\begin{align*}
			|J(U)-J(\tilde{U})|&=|J(U)-J(U_h^{(2)})+\eta|.
			\end{align*}
			Using the triangle inequality and inverse triangle inequality $|x+y|\geq|x|-|y|$ for $x,y \in \mathbb{R}$, we immediately get
			\begin{align*}
			|J(U)-J(\tilde{U})|&\leq|J(U)-J(U_h^{(2)})|+|\eta|\leq b_h|J(U)-J(\tilde{U})|+|\eta|,
			\end{align*}
			and 
			\begin{align*}
			|J(U)-J(\tilde{U})|&\geq-|J(U)-J(U_h^{(2)})|+|\eta|\geq -b_h|J(U)-J(\tilde{U})|+|\eta|.
			\end{align*}
			This implies that 
			$$(1-b_h)|J(U)-J(\tilde{U})| \leq |\eta| \leq (1+b_h)|J(U)-J(\tilde{U})|.$$
			This inequality is equivalent to 
			$$\frac{1}{1+b_h}|\eta|\leq|J(U)-J(\tilde{U})|\leq \frac{1}{1-b_h}|\eta|.$$
		\end{proof}
	\end{Corollary}
	As localization for the discretization error estimator $\eta_h$, we use the partition of unity (PU) technique proposed in \cite{RiWi15_dwr}.
	Here, we use the $Q_1$ finite elements as partition of unity. Let $\Phi_i \in Q_1$. Then, we know that $\sum_{i=1}^{|Q_1|}\Phi_i=1$.
	Of course, 
	\begin{align*}
	\frac{1}{2}\left(\eta_p(\tilde{U})(Z_h^{(2)}-\tilde{Z})+\eta_a(\tilde{U},\tilde{Z})(U_h^{(2)}-\tilde{U})\right),
	\end{align*}
	is equal to 
	\begin{align}
	&\frac{1}{2}\left(\eta_p(\tilde{U})((Z_h^{(2)}-\tilde{Z})\sum_{i=1}^{|Q_1|}\Phi_i)+\eta_a(\tilde{U},\tilde{Z})((U_h^{(2)}-\tilde{U})\sum_{i=1}^{|Q_1|}\Phi_i)\right) \nonumber\\
	=&\sum_{i=1}^{|Q_1|}\underbrace{\frac{1}{2}\left(\eta_p(\tilde{U})((Z_h^{(2)}-\tilde{Z})\Phi_i)+\eta_a(\tilde{U},\tilde{Z})((U_h^{(2)}-\tilde{U})\Phi_i)\right)}_{:=\eta_i^{PU}}. \label{eq: eta_iPU}
	\end{align}	 
	These $\eta_i^{PU}$ represent an error distribution of the partition of unity or their corresponding nodes respectively. This nodal error is distributed equally to all elements sharing this node.
	For more information about this, we refer to \cite{EnLaWi18,EndtLaWi20}.	 
	Finally, we want to show efficiency and reliability for $\eta_h$ using an additional assumption.
	\begin{Corollary}[see \cite{EndtLaRiSchafWi24_book_chapter}]
		Let us assume that all assumptions of Corollary~\ref{Corr: efficiency and reliablity eta} are fulfilled and that there exist a constant $\xi \in (0,1)$ such that
		$$|\eta_k+ \mathcal{R}_h^{(3),(2)}| \leq \xi |\eta_h|.$$
		Then the discretization error estimator $\eta_h$ is efficient and reliable with the constants $\overline{c}_h=\frac{(1+\xi)}{1-b_h}$ and $\underline{c}_h=\frac{(1-\xi)}{1+b_h}$.
		\begin{proof}
			This proof is also presented in \cite{EndtLaRiSchafWi24_book_chapter}.
			From the proof Corollary~\ref{Corr: efficiency and reliablity eta}, we know
			$$\frac{1}{1+b_h}|\eta|\leq|J(U)-J(\tilde{U})|\leq \frac{1}{1-b_h}|\eta|.$$
			Now we use that $|\eta|=|\eta_h+\eta_k+ \mathcal{R}_h^{(3),(2)}|$.
			From the triangle inequality and the inverse triangle inequality, we get
			\begin{align*}
			|\eta|&\leq |\eta_h|+|\eta_k+ \mathcal{R}_h^{(3),(2)}|\leq (1+\xi)|\eta_h|,\\
			|\eta|&\geq |\eta_h|-|\eta_k+ \mathcal{R}_h^{(3),(2)}|\geq (1-\xi)|\eta_h|.
			\end{align*}
			Combining this result with the result from above gives 
			$$\frac{(1-\xi)}{1+b_h}|\eta_h|\leq|J(U)-J(\tilde{U})|\leq \frac{(1+\xi)}{1-b_h}|\eta_h|.$$
			This already proves the result with the constants $\overline{c}_h=\frac{(1+\xi)}{1-b_h}$ and $\underline{c}_h=\frac{(1-\xi)}{1+b_h}$.
		\end{proof}
	\end{Corollary}
	\begin{Remark}
		Using some further estimates 
                \twick{from \cite{EndtLaRiSchafWi24_book_chapter}}
                one can show that the effectivity index 
		\begin{equation}\label{eq: def Ieff}
		I_{eff}:=\frac{|\eta_h|}{|J(U)-J(\tilde{U})|}
		\end{equation}
		satisfies $I_{eff}\in [\frac{(1-\xi)}{1+b_h},\frac{(1+\xi)}{1-b_h}]$.
		Additionally, we define the primal effectivity index  	
		\begin{equation*}\label{eq: def Ieffp}
		I_{eff,p}:=\frac{|\eta_p|}{|J(U)-J(\tilde{U})|},
		\end{equation*}
		and  the adjoint effectivity index 
		\begin{equation*}\label{eq: def Ieffa}
		I_{eff,a}:=\frac{|\eta_a|}{|J(U)-J(\tilde{U})|}.
		\end{equation*}
		We would like to mention that for linear problems $I_{eff,p}=I_{eff,a}$.
	\end{Remark}
	\begin{Remark}
		In our numerical experiments, we use Newton's method until $|\eta_k|<10^{-10}.$
	\end{Remark}

	\subsection{Multiple goals}
	In many applications, there are more functionals of interest. 
	For example in the case of our model problem, that could be the temperature at a given point $x_{\theta_1}$; i.e $J(U)=\theta(x_{\theta_1})$ or the mean magnitude of the velocity of the fluid; i.e $J(U)=\int_{\Omega} |v| \text{d}x$. Now let us assume we have $N \in \mathbb{N}$ functional evaluations $J_1, \ldots J_N$ in which we are interested in. One possible way would be to simply add all functionals to one, i.e $J_A:=\sum_{i=1}^{N}J_i$. However, this may lead to error cancellation effects. For instance, let us assume that $N=3$ and $J_1(U)-J_1(\tilde{U})=1, J_2(U)-J_2(\tilde{U})=1$ and $J_3(U)-J_3(\tilde{U})=-2$ then $J_A(U)-J_A(\tilde{U})=0$. This indicates that the error vanishes, but this is not true for any of the original functionals.

	\subsubsection{Combined Functional}
	To overcome the error cancellation effects when combining the functionals, we use a weighted sum with signed weights as done in \cite{HaHou03,Ha08,EndtWi17,EnLaNeWiWo20,EnLaWi18}.
	The combined functional, which is motivated by \cite{EnLaWi18} is finally given by
	\begin{equation} \label{eq: def Jc}
	J_c(v):=\sum_{i=1}^{N}w_iJ_i(v) \qquad \text{with} \qquad w_i:=\omega_i\frac{J(U_h^{(2)})-J(\tilde{U})}{|J(U_h^{(2)})-J(\tilde{U})|},
	\end{equation}
	where $w_i=\omega_i C$ with $C \in \mathbb{R}$ if $|J(U_h^{(2)})-J(\tilde{U})|=0$. Here, $w_i$ are user chosen weights. In our numerical examples, we set $w_i={|J(\tilde{U})|}^{-1}$ for non vanishing functionals and 1 for vanishing functionals. As explained in \cite{EnLaWi18} this choice of $w_i$ leads to a similar relative error in all the functionals. 
	For more details on this combined functional, we refer to \cite{EndtWi17,EnLaWi18}.
	Finally, we can apply the theory from the single case to the combined functional $J_c$.
	
	\subsection{Algorithms}
	In this subsection, we introduce our three main algorithms, namely 
	line search, the Newton's method, and the final adaptive algorithm.
	
	\begin{algorithm}[H]
\caption{Line search} \label{Alg: line search}
\begin{algorithmic}[1]	
	\Procedure{Line search}{$\beta, A, U, P$, maxiterL}      
	\State $\alpha\leftarrow 1$, $k\leftarrow 0$
	\While{$\Vert A(U+\alpha P)\Vert\geq\Vert A(U+\alpha P)\Vert$ \& $k$ $\geq$ maxiterL }
	\State $\alpha\leftarrow\beta \alpha$ 
	\State $k\leftarrow k+1$
	\EndWhile
	\Return $U+\alpha P$
	\EndProcedure
\end{algorithmic}	
\end{algorithm}
\begin{Remark}
	In our numerical experiments, we used $\beta=0.7$ and maxiterL=20.
\end{Remark}
\begin{algorithm}[H]
	\caption{Newton' method with line search}\label{Alg: Newton}
	\begin{algorithmic}[1]	
		\Procedure{Newton's Method}{$A, U^0$, res, maxiterN}     
		\State $k\leftarrow 0$, $U^k\leftarrow U^0$
		\While{$\Vert A(U^k)\Vert\geq$ res \& $k$ $\geq$ maxiterN }
		\State Solve $A'(U)(P)=-A(U)$ for $P$.
		\State $U^{k+1}$$\leftarrow$ Line search($\beta$, $A$, $U^k$, $P$, \text{maxiterL})
		\State $k\leftarrow k+1$
		\EndWhile
		\Return $U^k$
		\EndProcedure
	\end{algorithmic}	
\end{algorithm}
\begin{Remark}
	In our numerical experiments we used res=maximum($10^{-12}$,$10^{-9})\Vert A(U^0)\Vert$
\end{Remark}

\begin{algorithm}[H]
	\caption{The adaptive algorithm} \label{Alg: adapt. Alg}
	\begin{algorithmic}[1]
		\Procedure{Adaptive Algorithm}{$J$, $A$, $\mathcal{T}_0$, TOL, maxNDoFs}      
		\State $k\leftarrow 0$,$\mathcal{T}_k\leftarrow \mathcal{T}_0$, $\eta_h\leftarrow \infty$
		\While{$\eta_h$> $10^{-2}$ TOL \& $|\mathcal{T}_k|$ $\geq$ maxNDoFs }
		\State Solve \eqref{Goee: enriched model p} and \eqref{discretized weak form} to obtain $U_h^{(2)}$ and $U_h$ using Newton's method.
		\State Solve \eqref{Goee: enriched disc. adjoint} and \eqref{Goee: disc. adjoint} to obtain $Z_h^{(2)}$ and $Z_h$ using some linear solver.
		\State Compute $\eta_h$ and  the node-wise error contribution $\eta_i^{PU}$ as in \eqref{eq: eta_iPU}.
		\State Distribute $\eta_i^{PU}$ equally to all elements that share the node $i$.
		\State Mark the elements using some marking strategy.
		\State Refine the mesh
		\State $k\leftarrow k+1$
		\EndWhile
		\Return $J(U_h)$
		\EndProcedure
	\end{algorithmic}	
\end{algorithm}
\begin{Remark}
In Algorithm~\ref{Alg: adapt. Alg} and \ref{Alg: Newton} the arising linear problems are solved using the direct solver \newline UMFPACK~ \cite{UMFPACK}.
\end{Remark}
\begin{Remark}
In the marking strategy in Algorithm~\ref{Alg: adapt. Alg}, we mark to 10\% of all elements, when ordered with respect to their error contribution.
\end{Remark}
\begin{Remark}
If $J=J_c$ as defined in \eqref{eq: def Jc} then we return all $J_i(U_h)$.
\end{Remark}

	\section{Numerical experiments}
	\label{Sec: Num. Example}
	In this section, we conduct several numerical experiments in two space dimensions. The programming 
	code is based on the open source library deal.II \cite{deal2020,dealII91}
	and our algorithmic framework developed in \cite{EnLaWi18}.
	In all numerical experiments we use $\rho_0=998.21$,
	$\nu_0=2.216065960663198 \times 10^{-6}$ and 
	$E_A=14906.585117275014$.
	The thermal conductivity $k=0.5918$.
	Furthermore, we combine the single functionals to $J_c$ as defined in \eqref{eq: def Jc}.

	\subsection{Prerefinement around laser}
	\label{Pre-refinement for laser}
\begin{figure}[H]
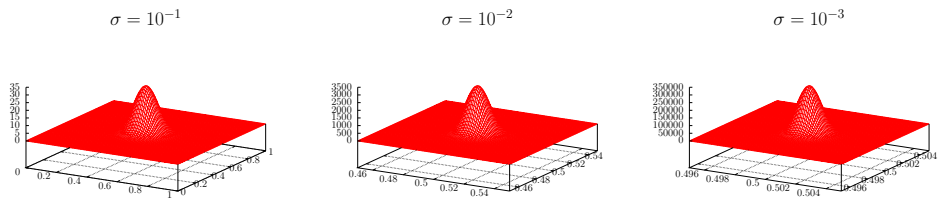

		{	\begin{minipage}{0.3\linewidth}
				\scalebox{0.35}{\input{Figures/Lasersource-12.tex}}
			\end{minipage}\hfill
			\begin{minipage}{0.3\linewidth}
				\scalebox{0.35}{\input{Figures/Lasersource-22.tex}}
			\end{minipage}\hfill	
			\begin{minipage}{0.3\linewidth}
				\scalebox{0.35}{\input{Figures/Lasersource-32.tex}}
			\end{minipage}
		}		
		\caption{$f_{E,y}(x)$ with $y=\frac{1}{2}(1,1)$ and $E=1$ on the domains $[0,1]^2$ for $\sigma=10^{-1}$ (left), $[0.45,0.55]^2$ for $\sigma=10^{-2}$ (center) and $[0.495,0.505]^2$ for $\sigma=10^{-3}$ (right). \label{fig: different sigmas}}
	\end{figure}

	In the numerical examples, we approximate a point laser in $y \in \mathbb{R}$ with the Gaussian function
	\begin{equation*}
	f_{E,y}(x):=E\left(\frac{1.0}{\sqrt{2\pi\sigma^2}} \right)^d e^{-\frac{(x-y)^2}{2\sigma^2}},
	\end{equation*}
	where $d$ is the space dimension, as energy source for the the heat equation.
	In the numerical experiments, we noticed that for large $E$ and small $\sigma$ prerefinement around $y$ is required. This is the case since the right hand side is not sufficiently resolved by the grid. In Figure~\ref{fig: different sigmas}, we show the right-hand side for $E=1$ and $\sigma \in \{10^{-1},10^{-2},10^{-3}\}$. Even though all sub figures look alike, we would like to mention that we zoom in each time we pick a smaller sigma. Furthermore, the maximum increases for smaller sigma.
	As a consequence also the solution of the heat equation for a given velocity field $v$ has high derivatives around $y$. The motivates 
	the choice of a mesh with local refinement around $y$. We would like to note that without prerefinement, either the right hand side was approximated badly on coarse grids, and on finer grids, Newton's method did not converge.

	\subsection{Temperature dependent thermal expansion coefficient}
	We construct the temperature dependent thermal expansion coefficient $\alpha(\theta)$ using linear splines with the data\footnote{\url{https://www.chemie.de/lexikon/Wasser_(Stoffdaten).html}} from Table~\ref{tab:termal expansion} as support points. Furthermore, we extrapolate $\alpha(\theta)$ in the exterior using the gradient at the boundaries.

	\begin{table}[htbp]
		\centering
		\begin{tabular}{|c|c|c|c|c|c|c|c|c|}
			\hline
			\textbf{temperature in $^\circ$C} & 0  & 4  & 5    & 10   & 15   & 20   & 25   & 30   \\
			\hline
			$\alpha=10^{-3}\times$             & -0.08 & 0 & 0.011 & 0.087 & 0.152 & 0.209 & 0.259 & 0.305 \\
			\hline
		\end{tabular}
		\begin{tabular}{|c|c|c|c|c|c|c|c|c|c|}
			\hline
			\textbf{temperature} & 35   & 40   & 45   & 50   & 60   & 70   & 80   & 90 & 99.63   \\
			\hline
			$\alpha=10^{-3}\times$                  & 0.347 & 0.386 & 0.423 & 0.457 & 0.522 & 0.583 & 0.64 & 0.696 &0.748 \\
			\hline
		\end{tabular}
		\caption{Thermal expansion coefficient for different temperatures.}
		\label{tab:termal expansion}
	\end{table}
	
	\subsection{Flow in a square with two lasers}	
	
	In the first example, there are two lasers in the square
	$\Omega:= (0, 0.3)^2$. We choose $\theta_0=293.15$ in \eqref{weak form}. The two lasers are represented with a source term, which is modeled by the sum of two Gaussian distributions from above. The centers of the lasers are at $x_1=(0.05,0.05)$ and $x_2=(0.25,0.05)$, which together with $\Omega$ is depicted in the left sub figure of Figure$~$\ref{fig: domain+mesh}. 
	The right hand side function is chosen as 
	\[
	f_E(x):=E\left( \frac{1}{\sqrt{2\pi\sigma^2}} \right)^2 e^{-\frac{(x-x_1)^2}{2\sigma^2}}+E\left( \frac{1}{\sqrt{2\pi\sigma^2}} \right)^2 e^{-\frac{(x-x_2)^2}{2\sigma^2}} 
	\]
	with $\sigma^{-2}=500$ and $E=10^4$.
	\subsubsection{Goal functionals}
	We are interested in the following goal  functionals:
	\begin{itemize}
		\item the mean value of the magnitude of the velocity:
		$$J_1(U) := \frac{1}{|\Omega|}\int_{\Omega}|v| \text{ d}x,$$
		\item the mean value of the temperature in the domain: 
		$$J_2(U) := \frac{1}{|\Omega|}\int_{\Omega}\theta dx,$$
		\item the temperature at  the center $y=(0.15,0.15)$ which is shown in Figure~\ref{fig: domain+mesh}:
		$$J_3(U) := \theta(y).$$
	\end{itemize}
	\subsubsection{{Discussion of our findings}}
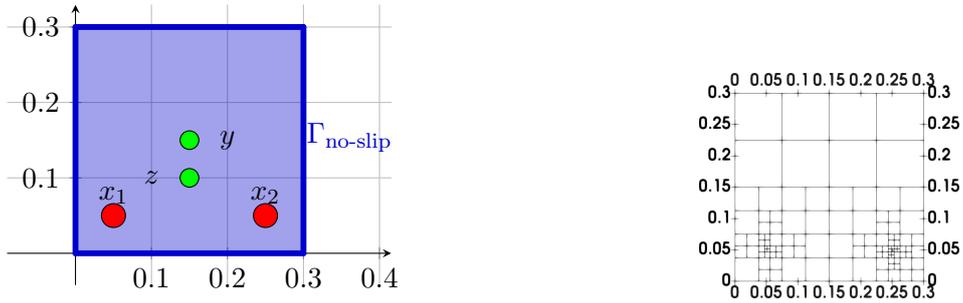
\begin{figure}[H]\label{fig: Qmega+Mesh}
			\definecolor{ffqqqq}{rgb}{1,0,0}
	\definecolor{qqqqcc}{rgb}{0,0,0.8}
	\begin{tikzpicture}[line cap=round,line join=round,>=triangle 45,x=10cm,y=10cm]
	\begin{axis}[
	x=10cm,y=10cm,
	axis lines=middle,
	ymajorgrids=true,
	xmajorgrids=true,
	xmin=-0.08951831248808961,
	xmax=0.4153884241802692,
	ymin=-0.04164612868197476,
	ymax=0.32934064961497966,
	xtick={-0.1,-0.0,...,0.4},
	ytick={-0.1,-0.0,...,0.4},]
	\clip(-0.08951831248808961,-0.04164612868197476) rectangle (0.4153884241802692,0.32934064961497966);
	\fill[line width=5.2pt,color=qqqqcc,fill=qqqqcc,fill opacity=0.37] (0,0) -- (0.3,0) -- (0.3,0.3) -- (0,0.3) -- cycle;
	\draw [line width=2.2pt,color=qqqqcc] (0,0)-- (0.3,0);
	\draw [line width=2.2pt,color=qqqqcc] (0.3,0)-- (0.3,0.3);
	\draw [line width=2.2pt,color=qqqqcc] (0.3,0.3)-- (0,0.3);
	\draw [line width=2.2pt,color=qqqqcc] (0,0.3)-- (0,0);
	\begin{scriptsize}
	\draw[color=qqqqcc] (0.3606880327483741,0.15339261920148556) node {$\Gamma_{\text{no-slip}}$};
	\draw [fill=green] (0.15,0.15) circle (3.5pt);
	\draw [color=black] (0.20,0.15) node {$y$};
	\draw [fill=green] (0.15,0.10) circle (3.5pt);
	\draw [color=black] (0.10,0.10) node {$z$};
	\draw [fill=ffqqqq] (0.05,0.05) circle (4.5pt);
	\draw[color=black] (0.05,0.077) node {$x_1$};
	\draw [fill=ffqqqq] (0.25,0.05) circle (4.5pt);
	\draw[color=black] (0.25,0.077) node {$x_2$};
	\end{scriptsize}
	\end{axis}
	\end{tikzpicture} \hfill
		\includegraphics[width=0.27\linewidth]{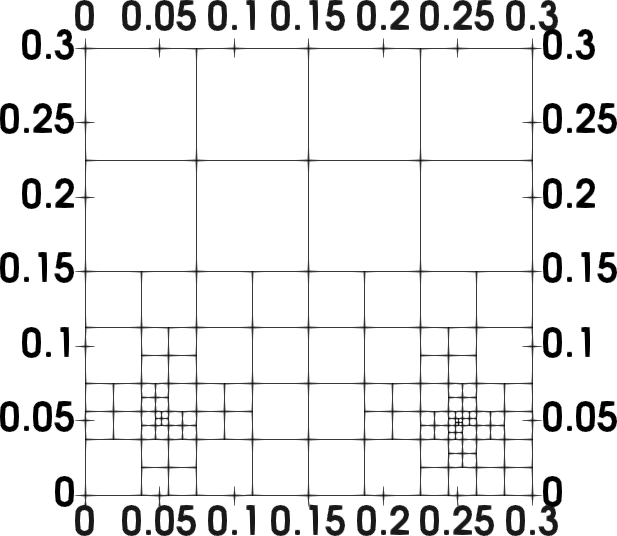} \hfill
		\caption{Example 1+2: The domain $\Omega$ $\Gamma_{\text{no-slip}}$ and $x_1$ and $x_2$ (left) and the initial mesh (right) \label{fig: domain+mesh}}
	\end{figure}

	As already mentioned in Subsection~\ref{Pre-refinement for laser} prerefinement at the position of the laser was required in order for the Newton's method to converge and to resolve the right the function $f_E$. The resulting initial mesh is visualized in Figure~\ref{fig: domain+mesh}. 
	\begin{figure}[H]
		\includegraphics[width=0.31\linewidth]{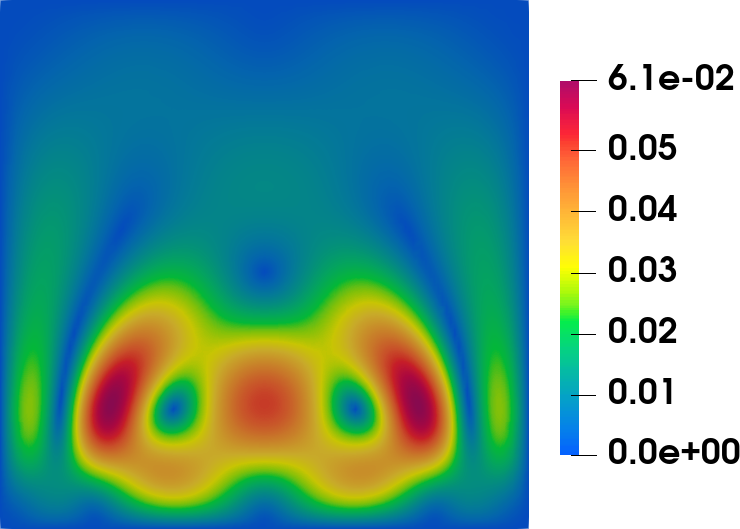} \hfill
		\includegraphics[width=0.31\linewidth]{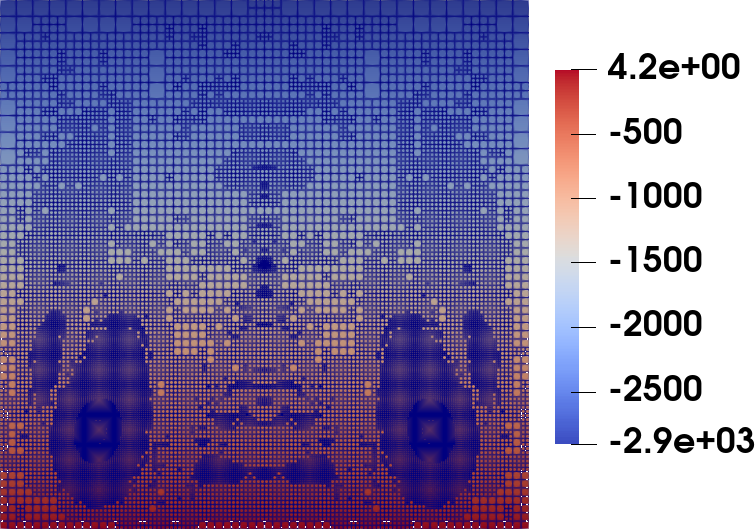} \hfill
		\includegraphics[width=0.31\linewidth]{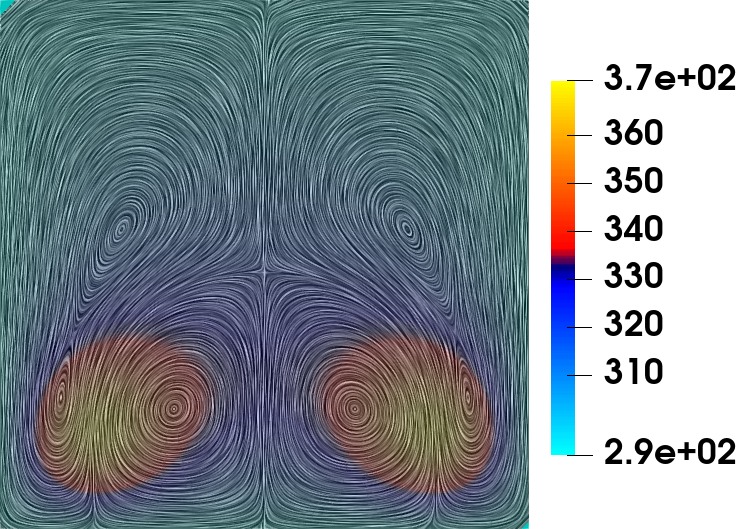}
		\caption{Example 1:  Magnitude of velocity(left), pressure and adaptively refined mesh (center) and streamlines and temperature (right)  \label{fig: Ex1: velo+flow + pressure+mesh+ temp}}
	\end{figure}
	{The numerical results for our first are presented in Figures~\ref{fig: Ex1: velo+flow + pressure+mesh+ temp}~-~\ref{fig: Ex1: errors} on pages~\pageref{fig: Ex1: velo+flow + pressure+mesh+ temp}~-~\pageref{fig: Ex1: errors}.}
	Both the streamlines (right) and the magnitude (left) of the velocity $v$ 
        are shown in Figure~\ref{fig: Ex1: velo+flow + pressure+mesh+ temp}. 
        Furthermore, there are at least six vortices. 

        However, refining in the corners of $\Omega$  generates more vortices at these corners. The impact of these vortices on our quantities of interest is negligible, since the goal oriented error estimator, would mark those elements for refinement otherwise.   
	In Figure \ref{fig: Ex1: errors}, we observe that the error in the combined functional is dominated by $J_1$. Therefore, the refinement is primarily driven by the refinement of this quantity of interest, which can be seen in Figure~\ref{fig: Ex1: velo+flow + pressure+mesh+ temp}. The relative errors of the other quantities of interest are almost one order of magnitude less. This also explains, why there is just a light refinement at $y$, compared to other examples with point evaluations as quantity of interest; see \cite{Endt21,RiWi15_dwr} and other examples in this work. The temperature $\theta$ in our computation is between $293.15K$ and $373.15K$ as represented in Figure~\ref{fig: Ex1: velo+flow + pressure+mesh+ temp}. We would like to mention that on coarse meshes, the effectivity index was between $0.1$ and $0.2$ on coarse meshes. This is an additional indicator the  the function $f_E$, was not good enough resolved on the initial mesh. This is shown in Figure~\ref{fig: Ex1: Ieff}. On finer meshes, we obtain effectivity indices between $0.68$ and $0.82$.
	\begin{figure}[H]

		{	\begin{minipage}{0.475\linewidth}
				\scalebox{0.49}{\input{Figures/IeffEx12.tex}}
				\caption{Example 1: Effectivity index. \label{fig: Ex1: Ieff}}
			\end{minipage}\hfill
			\begin{minipage}{0.475\linewidth}
				\scalebox{0.49}{\input{Figures/ErrorEx12.tex}}
				\caption{Example 1: Relative errors for $J_1$, $J_2$, $J_3$ and absolute error for $J_c$. \label{fig: Ex1: errors}}
			\end{minipage}	
		}	
		
	\end{figure}

	\subsection{Flow in a square with two lasers with different focus}
	In this subsection, we change the focus of the laser, namely the parameter $\sigma$. We use $$\sigma \in\{10^{-1},10^{-2},10^{-3},10^{-4}\}.$$ 
Furthermore in all examples, we set $E=100$. All other parameters are the same as in previous example.
	We have the following quantities of interest:
	\begin{itemize}
		\item the mean value of the magnitude of the velocity:
		$$J_1(U) := \frac{1}{|\Omega|}\int_{\Omega}|v| dx,$$
		\item the mean value of the temperature in the domain: 
		$$J_2(U) := \frac{1}{|\Omega|}\int_{\Omega}\theta dx,$$
		\item the temperature at $y$:
		$$J_3(U) := \theta(y),$$
		\item the temperature at $z$
		$$J_4(U) := \theta(z),$$
	\end{itemize}
	where $y=(0.15,0.15)$ and $z=(0.15,0.1)$ as visualized in Figure~\ref{fig: domain+mesh}, respectively.
	In our numerical examples we used the following reference values given Table~\ref{tab: Ex 2: reference values}.

	\begin{table}[H]
          \begin{center}
\scalebox{0.7}{		\begin{tabular}{|c|r|r|r|r|}
			\hline
			\textbf{Reference values} & { \centering$\sigma=10^{-1}$}        & { \centering$\sigma=10^{-2}$}    & { \centering$\sigma=10^{-3}$}       & { \centering$\sigma=10^{-4}$}      \\ \hline
			$J_1(U)=$                     & 4.568690$ \times10^{-3}$ & 4.821760$ \times10^{-3}$ & 5.692103$ \times10^{-3}$ & 5.705407$ \times10^{-3}$ \\ \hline
			$J_2(U)=$                     & 299.0619 & 302.1829 & 302.3692 & 302.3710 \\ \hline
			$J_3(U)=$                     & 305.7956 & 303.4749 & 303.4754 & 303.4754 \\ \hline
			$J_4(U)=$                     & 306.8565 & 307.3140 & 307.3146 & 307.3146 \\ \hline
		\end{tabular}}
		\caption{Example 2: reference values \label{tab: Ex 2: reference values}}
\end{center}
	\end{table}

\begin{figure}		
		\begin{figure}[H]
			\includegraphics[width=0.31\linewidth]{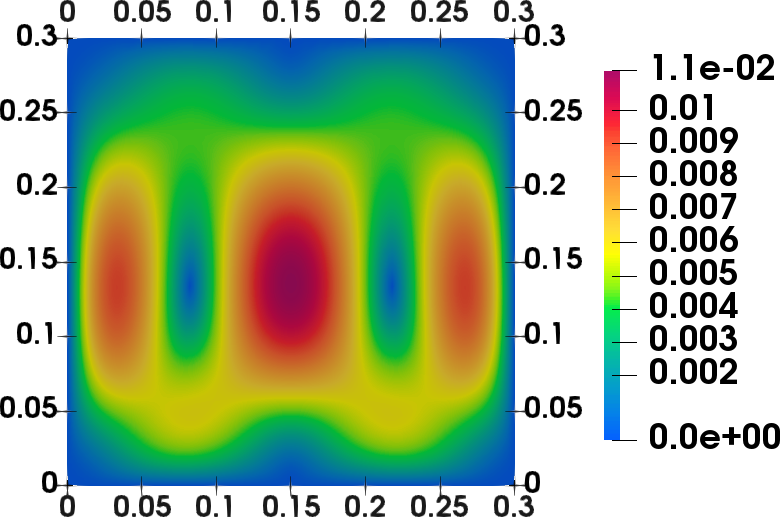}\hfill
			\includegraphics[width=0.31\linewidth]{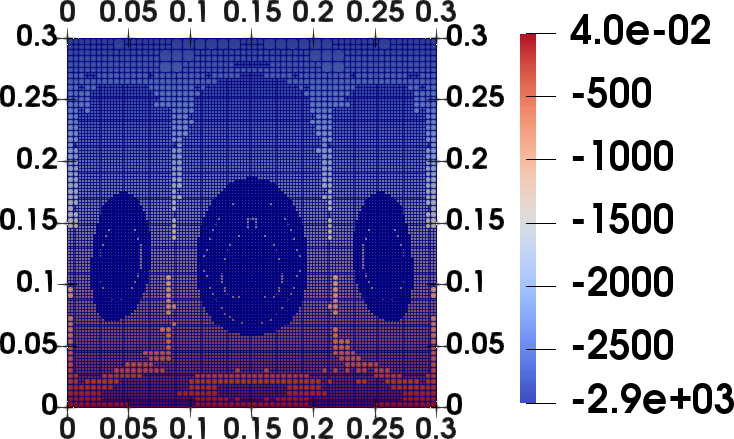}\hfill
			\includegraphics[width=0.31\linewidth]{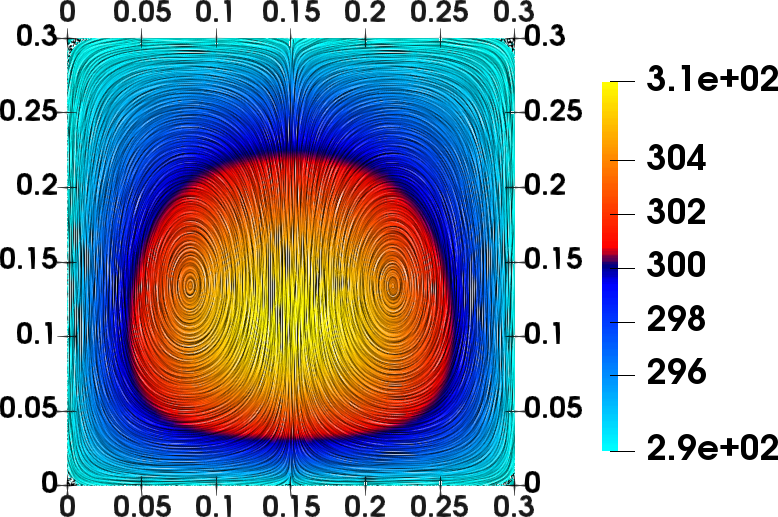} \\
			\caption{Example 2: Magnitude of velocity $v$ (left),  the pressure with the adaptive mesh(center) and temperature and flow field (right) for $\sigma=10^{-1}$. \label{fig: Ex2: sigma-1: Magvelo + Temp+flow} }
		\end{figure}
		\begin{figure}[H]
			\includegraphics[width=0.31\linewidth]{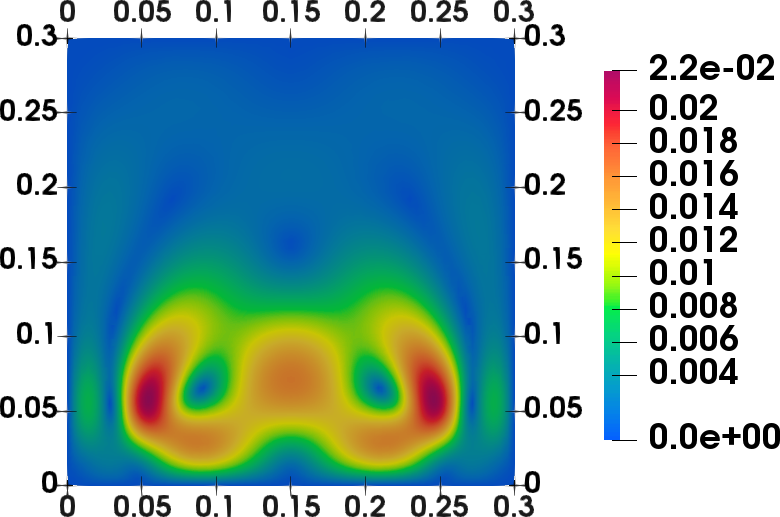}\hfill
			\includegraphics[width=0.31\linewidth]{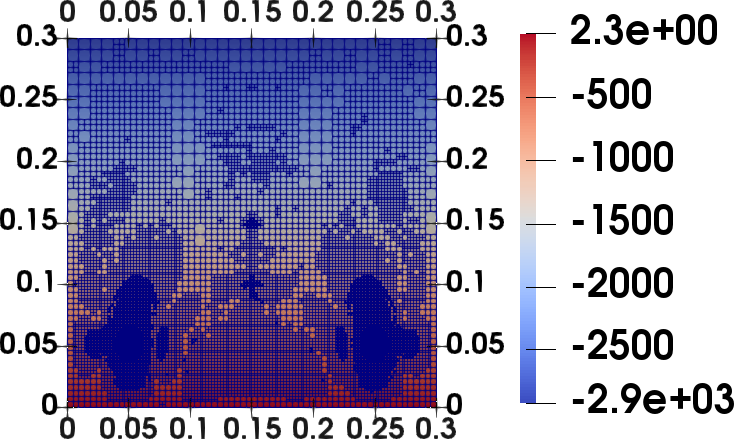}\hfill
			\includegraphics[width=0.31\linewidth]{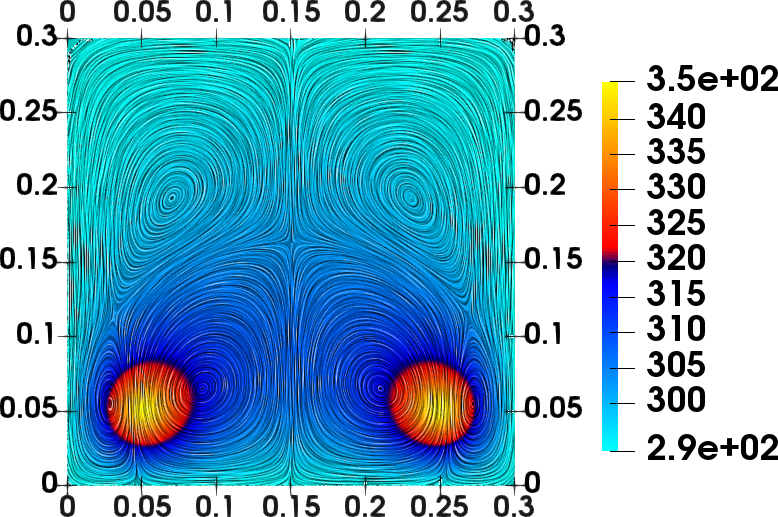} \\
			\caption{Example 2: Magnitude of velocity $v$ (left),  the pressure with the adaptive mesh(center) and temperature and flow field (right) for $\sigma=10^{-2}$. \label{fig: Ex2: sigma-2: Magvelo + Temp+flow}}
		\end{figure}
		\begin{figure}[H]
			\includegraphics[width=0.31\linewidth]{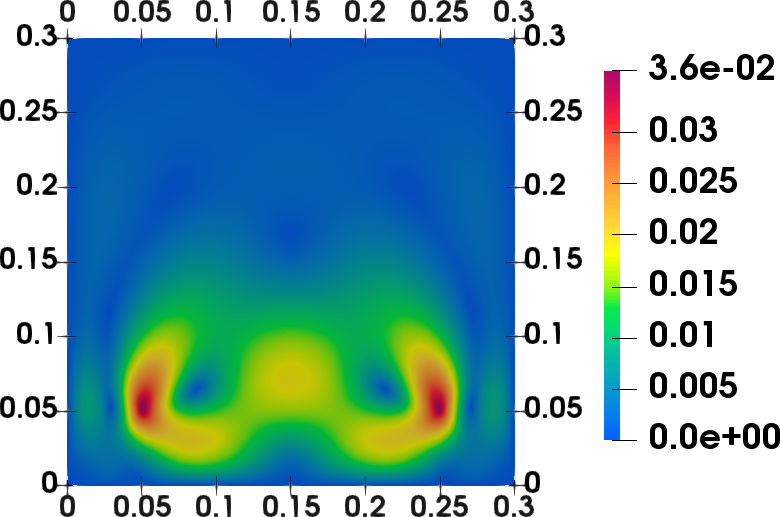}\hfill
			\includegraphics[width=0.31\linewidth]{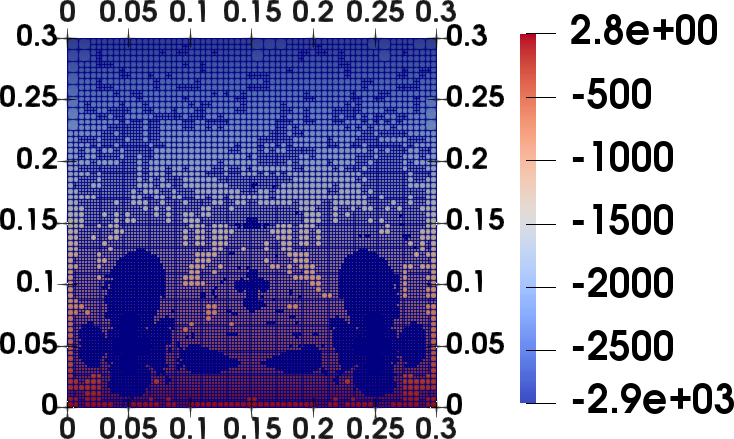}\hfill
			\includegraphics[width=0.31\linewidth]{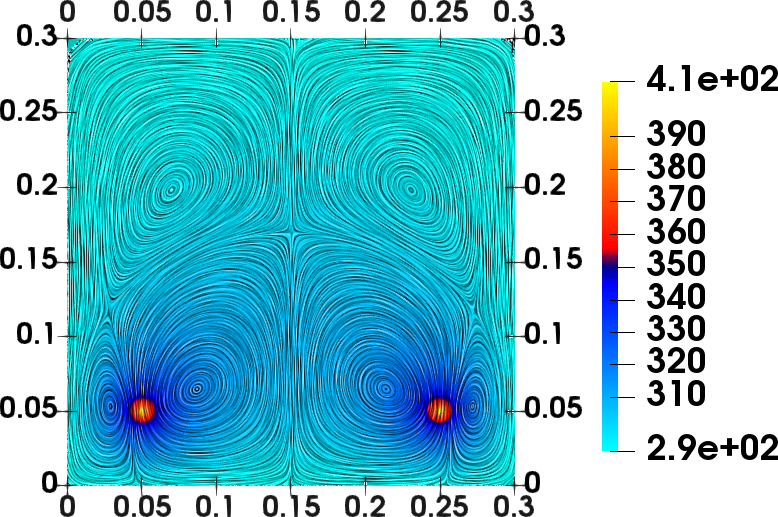} \\
			\caption{Example 2: Magnitude of velocity $v$ (left), the pressure with the adaptive mesh(center) and temperature and flow field (right) for $\sigma=10^{-3}$. \label{fig: Ex2: sigma-3: Magvelo + Temp+flow} }
		\end{figure}
		\begin{figure}[H]
			\includegraphics[width=0.31\linewidth]{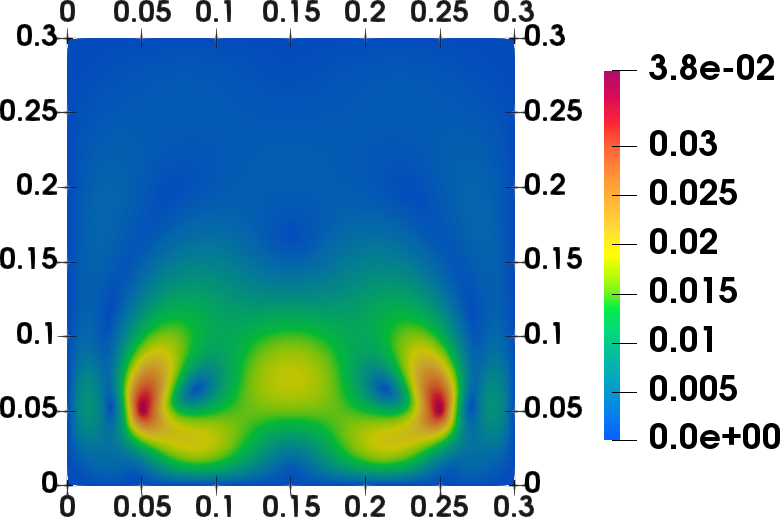}\hfill
			\includegraphics[width=0.31\linewidth]{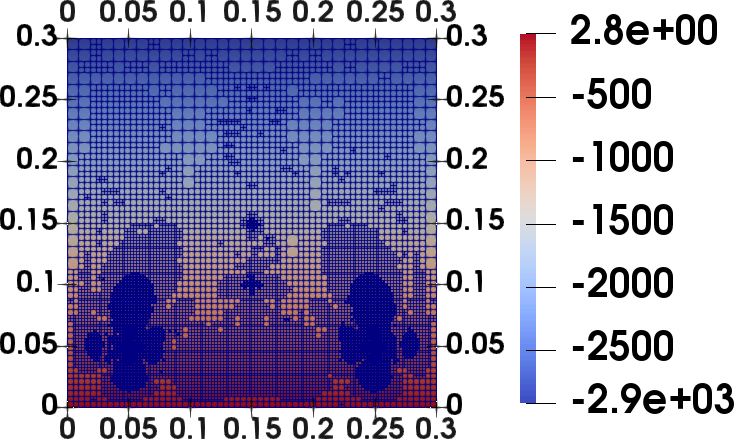}\hfill
			\includegraphics[width=0.31\linewidth]{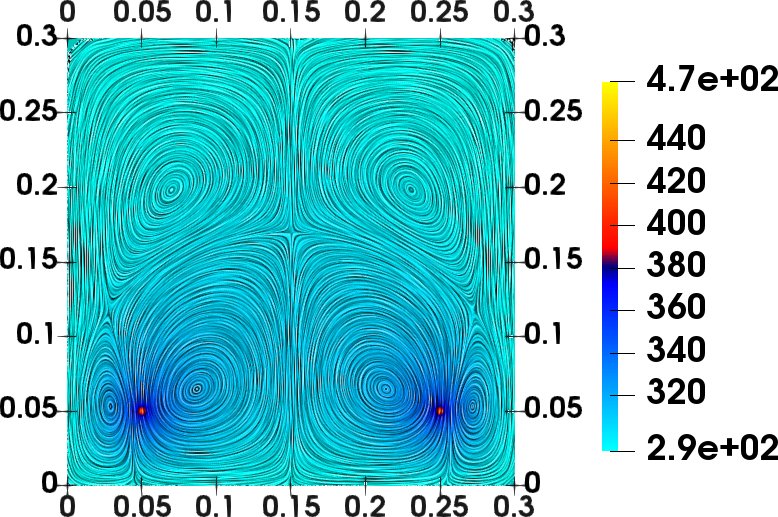} \\
			\caption{Example 2: Magnitude of velocity $v$ (left), the pressure with the adaptive mesh(center) and temperature and flow field (right) for $\sigma=10^{-4}$. \label{fig: Ex2: sigma-4: Magvelo + Temp+flow} }
		\end{figure}
	\end{figure}

	\bernhard{The numerical results for this example are presented in Figures~\ref{fig: Ex2: sigma-1: Magvelo + Temp+flow}~-~\ref{fig: Ex2: sigma-4: Errors} on pages~\pageref{fig: sigma-1 Ieff}~-~\pageref{fig: Ex2: sigma-4: Errors}.}
	For $\sigma=0.1$, the maximum of the temperature is close to the center, cf. Figure~\ref{fig: Ex2: sigma-1: Magvelo + Temp+flow}. This is in contrast to $\sigma \leq 0.01$, cf. Figures~\ref{fig: Ex2: sigma-2: Magvelo + Temp+flow}~-~~\ref{fig: Ex2: sigma-4: Magvelo + Temp+flow}, where the maximum is close to the point sources.
	This also changes the pattern of the flow, depicted on the right subfigures.
	Whereas for $\sigma=10^{-1}$, we just observe two vortices. 
	In Figures~\ref{fig: Ex2: sigma-2: Magvelo + Temp+flow} to \ref{fig: Ex2: sigma-4: Magvelo + Temp+flow}, we observe again six for $\sigma=10^{-2}$, $\sigma=10^{-3}$ and $\sigma=10^{-4}$. 
	For $\sigma=0.1$, cf. Figure~\ref{fig: Ex2: sigma-1: Magvelo + Temp+flow},  the magnitude of the velocity also has its maximum close to the center. 
	Moreover, the gradients of the temperature, and the largest temperature increase for smaller $\sigma$ as expected.
	For all $\sigma$, the functional $J_1$ dominates, cf. Figures~\ref{fig: sigma-1 Error}~-~~\ref{fig: Ex2: sigma-4: Errors}.
	Despite of that , the mesh for $\sigma=10^{-1}$ greatly differs from the others. The adaptively refined meshes of $\sigma=10^{-2}$, $\sigma=10^{-3}$ and $\sigma=10^{-4}$ share many similarities 
	like clover shaped refinement at $x_1$ and $x_2$. 
	Closer inspections also reveals this clover shaped refinement at $y$ and $z$, which is the typical refinement for point evaluations at this points, i.e. for $J_3$ and $J_4$. 
	This shows that even though $J_1$ is the dominating, our method provides refinement for the other functionals as well. 
	Furthermore, for the magnitude of velocity on the left of in Figures~\ref{fig: Ex2: sigma-1: Magvelo + Temp+flow}~-~\ref{fig: Ex2: sigma-4: Magvelo + Temp+flow} increases as $\sigma$ is decreasing.
	The effectivity indices are about 0.8 for $\sigma \le 0.01$, see Figures~\ref{fig: Ex2: sigma-2: Ieff}~-~\ref{fig: Ex2: sigma-4: Ieff}, if a certain refinement level is reached. 
	For $\sigma=10^{-4}$ it is very inaccurate for a small amount of DoFs. This is due to the strong point source. The result for $\sigma=0.1$ is quite different.

	\begin{figure}[H]

		{	\begin{minipage}{0.475\linewidth}
				\scalebox{0.49}{\input{Figures/IeffEx2+1e-12.tex}}
				\caption{Example 2: Effectivity Indices for $\sigma=10^{-1}$. \label{fig: sigma-1 Ieff}}
			\end{minipage}\hfill
			\begin{minipage}{0.475\linewidth}
				\scalebox{0.49}{\input{Figures/IeffEx2+1e-22.tex}}
				\caption{Example 2: Effectivity Indices for $\sigma=10^{-2}$.\label{fig: Ex2: sigma-2: Ieff}}
			\end{minipage}\newline
			\begin{minipage}{0.475\linewidth}
				\scalebox{0.49}{\input{Figures/IeffEx2+1e-32.tex}}
				\caption{Example 2: Effectivity Indices for $\sigma=10^{-3}$. \label{fig: Ex2: sigma-3: Ieff}}
			\end{minipage}\hfill
			\begin{minipage}{0.475\linewidth}
				\scalebox{0.49}{\input{Figures/IeffEx2+1e-42.tex}}
				\caption{Example 2: Effectivity Indices for $\sigma=10^{-4}$.\label{fig: Ex2: sigma-4: Ieff}}
			\end{minipage}
			
		}	
		
	\end{figure}

	\begin{figure}[H]

		{				
			\begin{minipage}{0.475\linewidth}
				\scalebox{0.49}{\input{Figures/ErrorsEx2+1e-12.tex}}
				\caption{Example 2: Relative errors in $J_{i}$, $i\in\{1,2,3,4\}$, and absolute error in $J_c$ for $\sigma=10^{-1}$. \label{fig: sigma-1 Error}}
			\end{minipage}	\hfill
			\begin{minipage}{0.475\linewidth}
				\scalebox{0.49}{\input{Figures/ErrorsEx2+1e-22.tex}}
				\caption{Example 2: Relative errors in $J_{i}$,  $i\in\{1,2,3,4\}$, and absolute error in $J_c$ for $\sigma=10^{-2}$. \label{fig: Ex2: sigma-2: Errors}}
			\end{minipage}
			\begin{minipage}{0.475\linewidth}
				\scalebox{0.49}{\input{Figures/ErrorsEx2+1e-32.tex}}
				\caption{Example 2: Relative errors in $J_{i}$,  $i\in\{1,2,3,4\}$, and absolute error in $J_c$ for $\sigma=10^{-3}$. \label{fig: Ex2: sigma-3: Errors}}
			\end{minipage}	\hfill
			\begin{minipage}{0.475\linewidth}
				\scalebox{0.49}{\input{Figures/ErrorsEx2+1e-42.tex}}
				\caption{Example 2: Relative errors in $J_{i}$,  $i\in\{1,2,3,4\}$, and absolute error in $J_c$ for $\sigma=10^{-4}$. \label{fig: Ex2: sigma-4: Errors}}
			\end{minipage}
		}	
		
	\end{figure}

	\subsection{Circle with a comparison between density and Boussinesq}
	
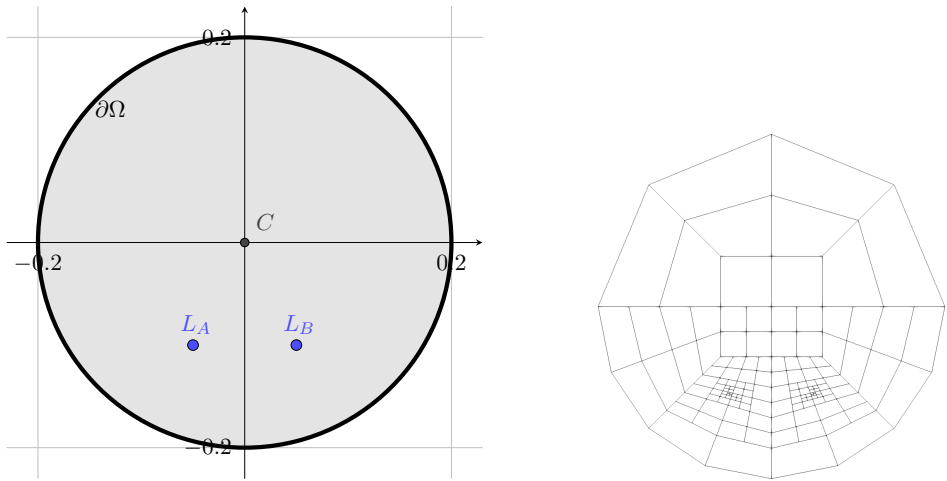
\begin{figure}[htbp]
		
	{	\scalebox{0.8}{\definecolor{ududff}{rgb}{0.30196078431372547,0.30196078431372547,1}
\definecolor{uuuuuu}{rgb}{0.26666666666666666,0.26666666666666666,0.26666666666666666}
\definecolor{uuuuu2u}{rgb}{0.06666666666666666,0.06666666666666666,0.06666666666666666}
\begin{tikzpicture}[line cap=round,line join=round,>=triangle 45,x=1cm,y=1cm]
\begin{axis}[
x=17cm,y=17cm,
axis lines=middle,
ymajorgrids=true,
xmajorgrids=true,
xmin=-0.23,
xmax=0.23,
ymin=-0.23,
ymax=0.23,
xtick={-0.2,-0.0,...,0.2},
ytick={-0.2,-0.0,...,0.0,0.2},]
\clip(-0.23,-0.23) rectangle (0.23,0.23);
\draw [line width=2pt] (0,0) circle (0.2);
\draw [fill=uuuuu2u,fill opacity=0.11] (0,0) circle (0.2);
\begin{scriptsize}
\draw [fill=uuuuuu] (0,0) circle (2pt);
\draw[color=uuuuuu] (0.02,0.02) node {$C$};
\draw[color=black] (-0.13,0.13) node {$\partial \Omega$};
\draw [fill=ududff] (-0.05,-0.1) circle (2.5pt);
\draw[color=ududff] (-0.047258124372038975,-0.08) node {$L_A$};
\draw [fill=ududff] (0.05,-0.1) circle (2.5pt);
\draw[color=ududff] (0.05266718819529195,-0.08) node {$L_B$};
\end{scriptsize}
\end{axis}
\end{tikzpicture} }\hfill
		\includegraphics[width=0.36\linewidth]{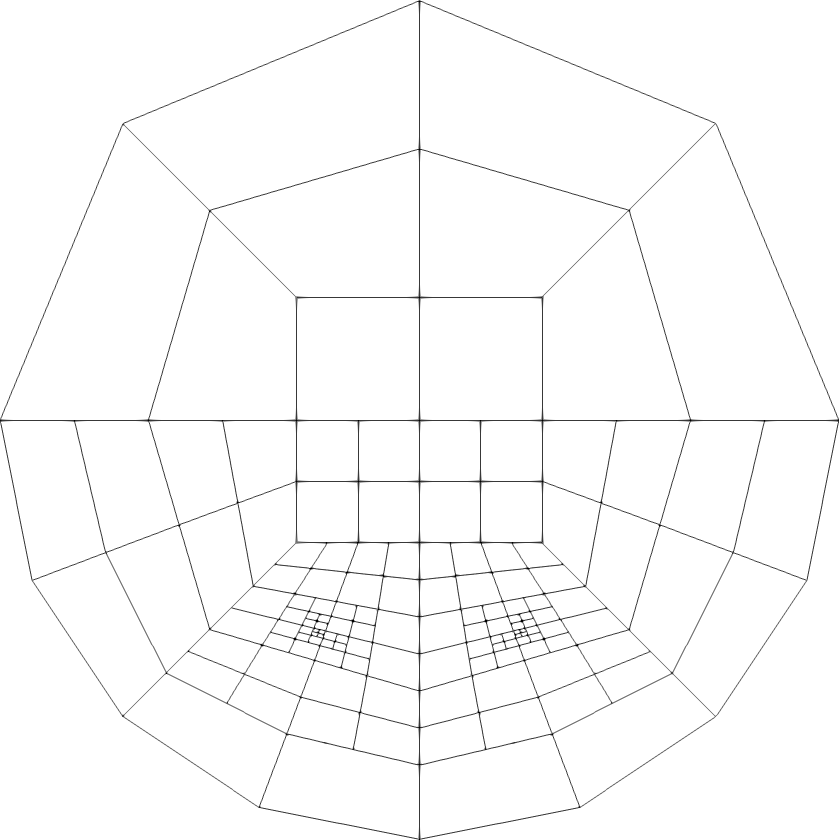}
		\caption{The domain $\Omega$, the center $C$ and the source points of the two lasers $L_A$ and $L_B$ on the left and the initial mesh with prerefinement on $L_A$ and $L_B$ on the right. \label{Fig:Ex:Circle:Omega+InitMesh}}}
	\end{figure}

	In the last example, the computational domain $\Omega$ is a circle with radius $0.2$ and center $C=(0,0)$, i.e. $\Omega:=\{x \in \mathbb{R}^2: |x|^2 < 0.04\}$. Additionally, we have two different lasers on the points $L_A:=\frac{1}{20}(-1,-2)$ and $L_B:=\frac{1}{20}(1,-2)$ with focus $\sigma=10^{-3}$.
	These quantities are depicted on Figure~\ref{Fig:Ex:Circle:Omega+InitMesh}. 
	We are interested in the following functionals
	\begin{align*}
	J_1(U)=v_1(L_A), \qquad & J_2(U)=v_2(L_A), \\
	J_3(U)=|v|^2(L_A), \qquad & J_4(U)=p(L_A)-p(L_B), \\
	J_5(U)=\theta(C), \qquad & J_6(U)=\int_{\partial\Omega}\nabla \theta \cdot n \text{ds}_x,
	\end{align*}
	where $n$ denotes the outer normal vector.
	Furthermore, we compare two different models in this example, \twick{which is 
        on the one hand our model problem~\eqref{eq:stationary_navierstokes}. In this example, we use the terminology 'Density'}, if we refer to our model problem. In case of the following  Boussinesq problem, we use the terminology 'Boussinesq' to refer to this problem.
	
\begin{figure}[htbp]
		\includegraphics[width=0.36\linewidth]{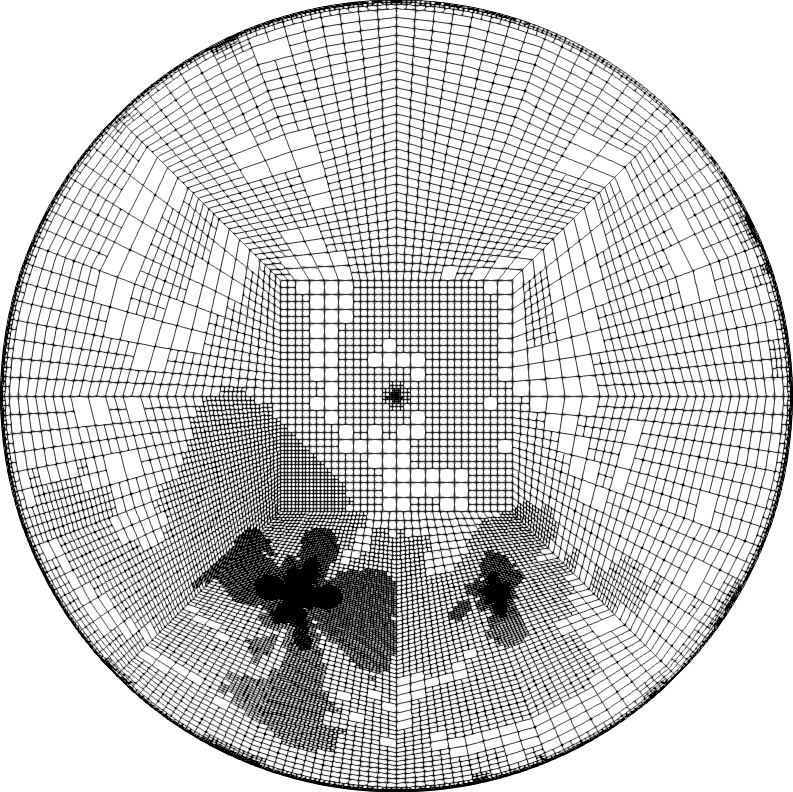}\hfill
		\includegraphics[width=0.36\linewidth]{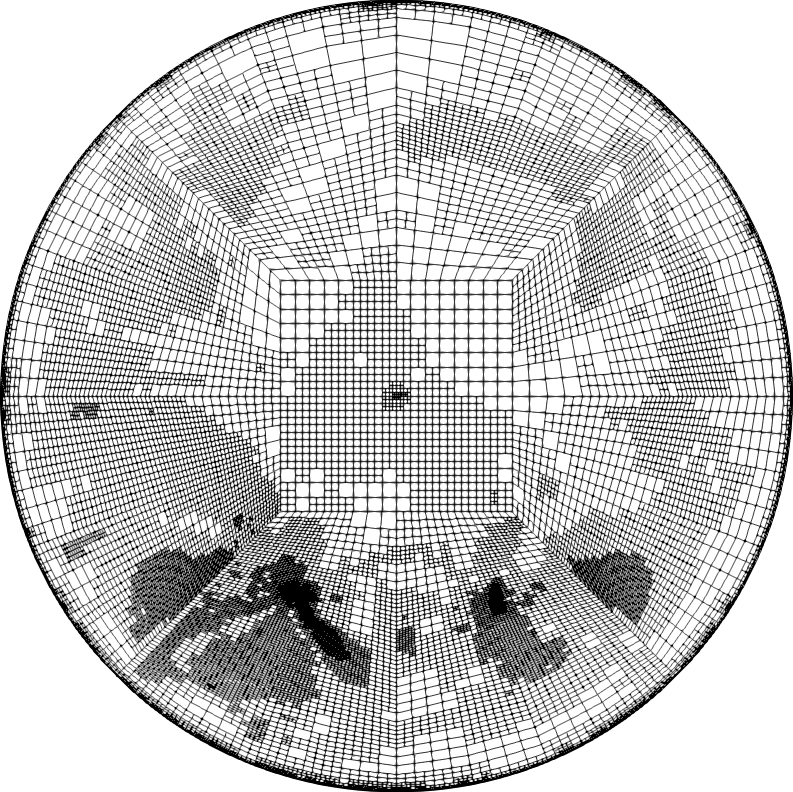} \\
		\caption{Example 3: Adaptively refined meshes for Density(left) and Boussinesq (right). \label{fig: Ex3: AdaptMesh}}
	\end{figure}

\begin{figure}[htbp]
		\includegraphics[width=0.39\linewidth]{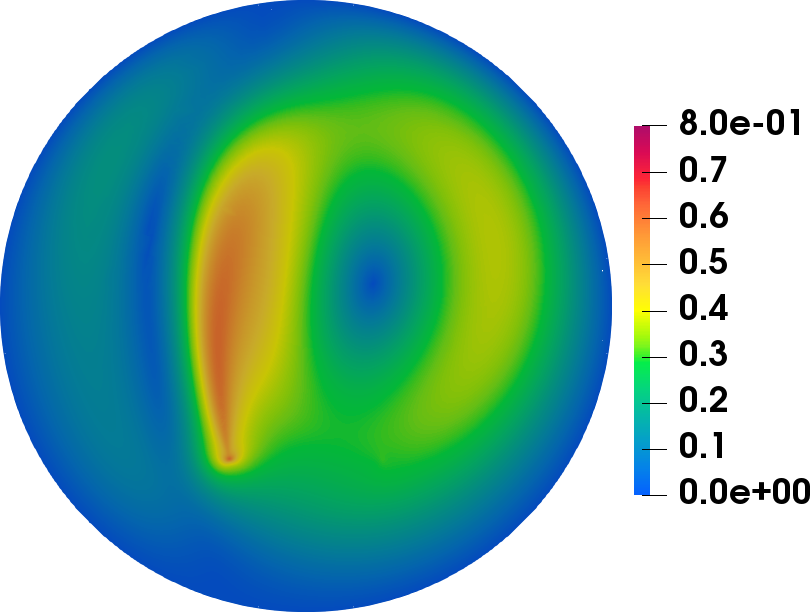}\hfill
		\includegraphics[width=0.39\linewidth]{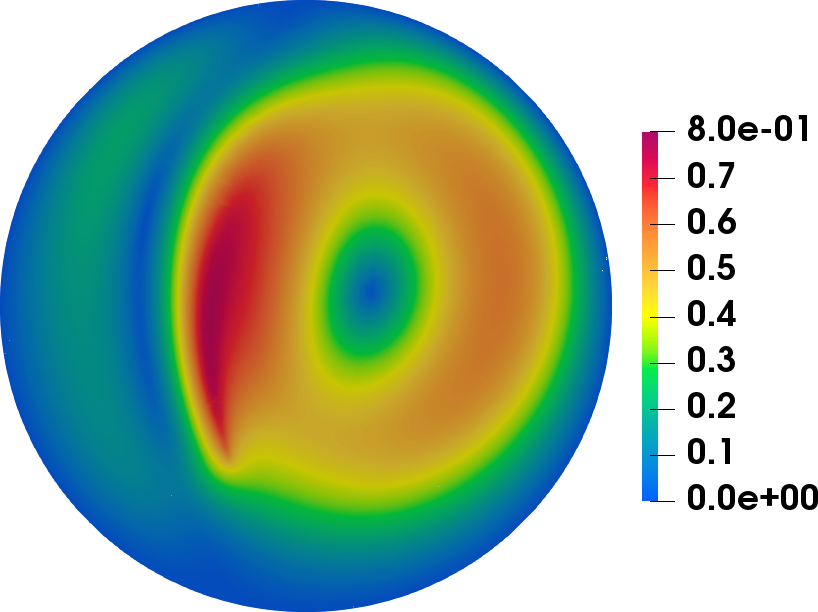} \\	
		\caption{Example 3: Magnitude of the velocity for Density(left) and Boussinesq (right). \label{fig: Ex3: MagV}}
	\end{figure}

\begin{figure}[htbp]
		\includegraphics[width=0.36\linewidth]{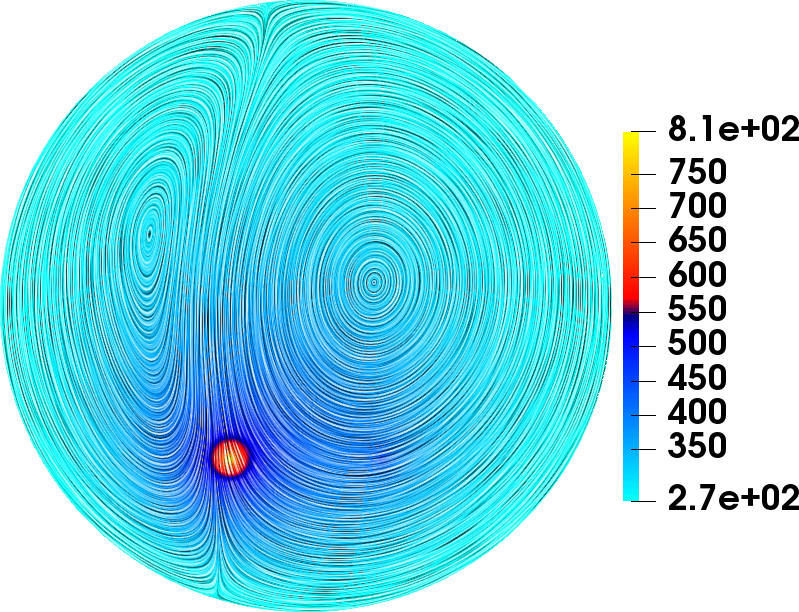}\hfill
		\includegraphics[width=0.36\linewidth]{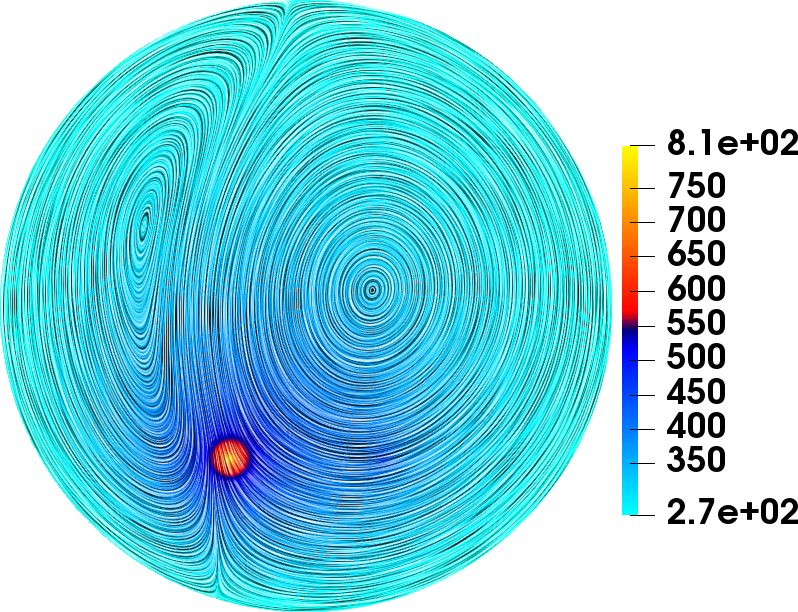} \\	
		\caption{Example 3: Streamlines of the flow and temperature for Density(left) and Boussinesq (right). \label{fig: Ex3: FlowT}}
	\end{figure}

	The Boussinesq problem {is a modification of \cite{LoBo96,BeuEndtLaWi23}} and given by: Find the velocity field $v$, the pressure $p$ and the temperature $\theta$ such that
	\begin{equation}
	\begin{aligned}
	-div\left(\nu(\theta)\rho(\nabla v + \nabla v^T\right) + (\rho v \cdot \nabla) v  - \nabla p - g \rho \alpha(\theta)(\theta-293.15) &=0,  \nonumber \\
	-div(v)&=0, \nonumber\\
	-div(k \nabla \theta) + v\cdot \nabla \theta & = f_{E},\label{eq:stationary_navierstokes_Boussinesq}\\
	\nu(\theta) &= \nu_0 e^{\frac{E_A}{R \theta}}.\nonumber
	\end{aligned}
	\end{equation}
	The boundary conditions for the Density model and Boussinesq are
	\begin{equation}
	\begin{aligned}
	v&=0  \qquad \text{ on }\partial \Omega, \\
	\theta&=274.15 \qquad  \text{ on }\partial \Omega.
	\end{aligned}
	\end{equation}
	Furthermore,  for both problems, i.e. the Density model and Boussinesq,
	\begin{equation*}
	f_E:=\gamma f_{E,L_A}+\gamma^{-1} f_{E,L_B},
	\end{equation*}
	with $E=200$ and $\gamma=2$. This choice of $\gamma$ makes the laser pointing at $L_A$ stronger then the one pointing on $L_B$, \Sven{see the left picture of Figure~\ref{Fig:Ex:Circle:Omega+InitMesh}.}
	The reference values used in the numerical experiments are given in Table~\ref{tab: Ex3: reference values}.
	\bernhard{The numerical results are presented in Figures~\ref{fig: Ex3: AdaptMesh}~-~\ref{fig: Ex3:Boussinesq :Comparison J4} on pages~\pageref{fig: Ex3: AdaptMesh}~-~\pageref{fig: Ex3:Boussinesq :Comparison J4} and Table~\ref{tab: Newtondata} on page~\pageref{tab: Newtondata}.}
	\Sven{As in the previous examples, the initial mesh has some prerefinement around the point sources, see the left picture of Figure~\ref{Fig:Ex:Circle:Omega+InitMesh}.
		Figure~\ref{fig: Ex3: FlowT} depicts the temperature and the streamlines of the flow. There are only minimal differences in the temperature.
		However, the magnitude of the velocity is larger for the  Boussinesq model, see Figure~\ref{fig: Ex3: MagV}.
		In both cases, there are two vortices, one on the left and one on the right. Whereas the right one looks quite similar in both cases, 
		the right one, see Figure~\ref{fig: Ex3: FlowT}, becomes more anisotropic for the  Boussinesq model.
		Also the adaptively refined meshes which are displayed in Figure~\ref{fig: Ex3: AdaptMesh} refine more from the point sources away for the  Boussinesq model.
		Figures~\ref{fig: Ex3: Errors Density} and \ref{fig: Ex3: Errors Boss} display the errors of the functionals $J_i$, $i=1,\ldots,6$. 
		The function $J_6$ dominates the error in both cases. In comparison to the previous examples, the dominance is not so strong.
		For coarser meshes, $J_4$ and in the case of  Boussinesq also $J_1$ are quite large. This is changing with the refinement process.
		For Boussinesq, then also $J_3$ will become quite large. This results in in very slow decrease of the functional $J_6$ in the first refinement steps which is worse than a uniform refinement, cf.
		Figures ~\ref{fig: Ex3: Density: Comparison J1} and \ref{fig: Ex3: Boussinesq: Comparison J1}, in particular for the density model. \bernhard{However, these figures also show that the adaptive procedure converges with approximately $O(\text{DoFs}^{-1})$, in comparison to $O(\text{DoFs}^{-\frac{1}{2}})$ for uniform refinement.}
		Figures ~\ref{fig: Ex3:Density :Comparison J4} and \ref{fig: Ex3:Boussinesq :Comparison J4} display the situations for the functional $J_4$.
		This functional is decreased quite rapidly in both cases.
		\bernhard{Here, both adaptive and uniform refinement have an error reduction in the order of $O(\text{DoFs}^{-1})$.}
		The effectivity indices, which are presented in Figures~\ref{fig: Ex3: Ieff Density} and ~\ref{fig: Ex3: Ieff Boss} are quite different for both models.
		The authors observe a effectivity indices greater than one for the Density model. They are about $1.2$.
		The effectivity indices are more away from one for the Boussinesq. Moreover, they  are less than one.
	}

	\begin{table}[htbp]
		\centering
		\begin{tabular}{|c|r|r|}
			\hline
			reference values& \multicolumn{1}{c|}{Density}  & \multicolumn{1}{c|}{Boussinesq}  \\ \hline
			$J_1$ & -0.12361494 & -0.16786688 \\ \hline
			$J_2$ &  0.57190820  & 0.48240847  \\ \hline
			$J_3$ &  0.34235964 & 0.26089722  \\ \hline
			$J_4$ &  8.62163631  & -73.2536833 \\ \hline
			$J_5$ & 353.526189  & 353.919688  \\ \hline
			$J_6$ & -844.882355 & -844.882273 \\ \hline
		\end{tabular}
		\caption{Example 3: Reference values \label{tab: Ex3: reference values}}
	\end{table}
	
		\begin{table}[htbp]
	\scalebox{0.75}{	\begin{tabular}{|r|rrr|rrr|}
			\hline
			\textbf{} & \multicolumn{3}{c|}{Density}                                                                               & \multicolumn{3}{c|}{Boussinesq}                                                                               \\ \hline
			\# refinements        & \multicolumn{1}{r|}{DoFs} & \multicolumn{1}{r|}{Newton steps} & {\# Linesearch} & \multicolumn{1}{r|}{DoFs} & \multicolumn{1}{r|}{Newton steps} & \multicolumn{1}{r|}{\# Linesearch}\\ \hline
			0         & \multicolumn{1}{r|}{1 680}         & \multicolumn{1}{r|}{15}            & 61                         & \multicolumn{1}{r|}{1 680}         & \multicolumn{1}{r|}{19}            & 88                         \\ \hline
			1         & \multicolumn{1}{r|}{2 458}         & \multicolumn{1}{r|}{4}             & 0                          & \multicolumn{1}{r|}{2 420}         & \multicolumn{1}{r|}{6}             & 0                          \\ \hline
			2         & \multicolumn{1}{r|}{3 216}         & \multicolumn{1}{r|}{4}             & 0                          & \multicolumn{1}{r|}{3 114}         & \multicolumn{1}{r|}{7}             & 0                          \\ \hline
			3         & \multicolumn{1}{r|}{4 428}         & \multicolumn{1}{r|}{4}             & 0                          & \multicolumn{1}{r|}{4 066}         & \multicolumn{1}{r|}{7}             & 0                          \\ \hline
			4         & \multicolumn{1}{r|}{5 982}         & \multicolumn{1}{r|}{3}             & 0                          & \multicolumn{1}{r|}{5 224}         & \multicolumn{1}{r|}{7}             & 0                          \\ \hline
			5         & \multicolumn{1}{r|}{7 940}         & \multicolumn{1}{r|}{4}             & 0                          & \multicolumn{1}{r|}{6 936}         & \multicolumn{1}{r|}{6}             & 0                          \\ \hline
			6         & \multicolumn{1}{r|}{10 404}        & \multicolumn{1}{r|}{3}             & 0                          & \multicolumn{1}{r|}{9 312}         & \multicolumn{1}{r|}{6}             & 0                          \\ \hline
			7         & \multicolumn{1}{r|}{13 532}        & \multicolumn{1}{r|}{3}             & 0                          & \multicolumn{1}{r|}{12 240}        & \multicolumn{1}{r|}{7}             & 0                          \\ \hline
			8        & \multicolumn{1}{r|}{17 636}        & \multicolumn{1}{r|}{3}             & 0                          & \multicolumn{1}{r|}{15 932}        & \multicolumn{1}{r|}{6}             & 0                          \\ \hline
			9        & \multicolumn{1}{r|}{23 032}        & \multicolumn{1}{r|}{3}             & 0                          & \multicolumn{1}{r|}{21 132}        & \multicolumn{1}{r|}{5}             & 0                          \\ \hline
			10        & \multicolumn{1}{r|}{29 862}        & \multicolumn{1}{r|}{3}             & 0                          & \multicolumn{1}{r|}{27 958}        & \multicolumn{1}{r|}{6}             & 0                          \\ \hline
			11        & \multicolumn{1}{r|}{38 774}        & \multicolumn{1}{r|}{3}             & 0                          & \multicolumn{1}{r|}{36 486}        & \multicolumn{1}{r|}{5}             & 0                          \\ \hline
			12        & \multicolumn{1}{r|}{50 320}        & \multicolumn{1}{r|}{3}             & 0                          & \multicolumn{1}{r|}{47 652}        & \multicolumn{1}{r|}{6}             & 0                          \\ \hline
			13        & \multicolumn{1}{r|}{65 480}        & \multicolumn{1}{r|}{3}             & 0                          & \multicolumn{1}{r|}{61 512}        & \multicolumn{1}{r|}{5}             & 0                          \\ \hline
			14        & \multicolumn{1}{r|}{84 822}        & \multicolumn{1}{r|}{3}             & 0                          & \multicolumn{1}{r|}{80 472}        & \multicolumn{1}{r|}{5}             & 0                          \\ \hline
			15        & \multicolumn{1}{r|}{110 122}       & \multicolumn{1}{r|}{3}             & 0                          & \multicolumn{1}{r|}{105 170}       & \multicolumn{1}{r|}{5}             & 0                          \\ \hline
			16        & \multicolumn{1}{r|}{143 414}       & \multicolumn{1}{r|}{3}             & 0                          & \multicolumn{1}{r|}{137 784}       & \multicolumn{1}{r|}{5}             & 0                          \\ \hline
			17        & \multicolumn{1}{r|}{186 514}       & \multicolumn{1}{r|}{3}             & 0                          & \multicolumn{1}{r|}{182 182}       & \multicolumn{1}{r|}{5}             & 0                          \\ \hline
			18        & \multicolumn{1}{r|}{243 002}       & \multicolumn{1}{r|}{3}             & 0                          & \multicolumn{1}{r|}{236 512}       & \multicolumn{1}{r|}{4}             & 0                          \\ \hline
			19        & \multicolumn{1}{r|}{315 344}       & \multicolumn{1}{r|}{3}             & 0                          & \multicolumn{1}{r|}{309 694}       & \multicolumn{1}{r|}{5}             & 0                          \\ \hline
			20        & \multicolumn{1}{r|}{409 738}       & \multicolumn{1}{r|}{2}             & 0                          & \multicolumn{1}{r|}{406 702}       & \multicolumn{1}{r|}{5}             & 0                          \\ \hline
		\end{tabular}}
		\caption{Comparison of required number of Newton and line search steps
			for the two different models: \# refinements denotes the current refinement level,DoFs denotes the number of degrees of freedom, Newton steps the number of Newton steps, \#~Linesearch denotes the total number of line searches steps required on the current refinement level \label{tab: Newtondata}}
	\end{table}

\newpage	
\mbox{}


	\begin{figure}[htbp]

		{	\begin{minipage}{0.4\linewidth}
				\scalebox{0.49}{\input{Figures/IeffExCircle_Density2.tex}}
				\caption{Example 3: Effectivity Indices for Density. \label{fig: Ex3: Ieff Density}}
			\end{minipage}\hfill
			\begin{minipage}{0.4\linewidth}
				\scalebox{0.49}{\input{Figures/IeffExCircle_Boussinesq2.tex}}
				\caption{Example 3: Effectivity Indices for Boussinesq. \label{fig: Ex3: Ieff Boss}}
			\end{minipage}	
		}	
		\end{figure}

\begin{figure}[htbp]

		{	\begin{minipage}{0.46\linewidth}
				\scalebox{0.49}{\input{Figures/ErrorExCircle_Density2.tex}}
				\caption{Example 3: Relative errors in $J_i$, $i=1,\ldots,6,$ and absolute error in $J_c$ for Density. \label{fig: Ex3: Errors Density}}
			\end{minipage}\hfill
			\begin{minipage}{0.46\linewidth}
				\scalebox{0.49}{\input{Figures/ErrorExCircle_Boussinesq2.tex}}
				\caption{Example 3: Relative errors in $J_i$, $i=1,\ldots,6,$ and absolute error in $J_c$ for Boussinesq. \label{fig: Ex3: Errors Boss}}
			\end{minipage}	
		}	

		{	\begin{minipage}{0.46\linewidth}
				\scalebox{0.49}{\input{Figures/Comp-J6-density2.tex}}
				\caption{Example 3: Relative errors in $J_{6}$ using adaptive and uniform(global) refinement for Density. \label{fig: Ex3: Density: Comparison J1}}
			\end{minipage}\hfill
			\begin{minipage}{0.46\linewidth}
				\scalebox{0.49}{\input{Figures/Comp-J6-boussi2.tex}}
				\caption{Example 3: Relative errors in $J_{6}$ using adaptive and uniform(global) refinement for Boussinesq. \label{fig: Ex3: Boussinesq: Comparison J1}}
			\end{minipage}	
		}	

		{	\begin{minipage}{0.46\linewidth}
				\scalebox{0.49}{\input{Figures/Comp-J4-density2.tex}}
				\caption{Example 3: Relative errors in $J_{4}$ using adaptive and uniform(global) refinement for Density. \label{fig: Ex3:Density :Comparison J4}}
			\end{minipage}\hfill
			\begin{minipage}{0.46\linewidth}
				\scalebox{0.49}{\input{Figures/Comp-J4-boussi2.tex}}
				\caption{Example 3: Relative errors in $J_{4}$ using adaptive and uniform(global) refinement for Boussinesq. \label{fig: Ex3:Boussinesq :Comparison J4}}
			\end{minipage}	
		}	
	\end{figure}
	
	\Sven{Both models are highly nonlinear and solved with Newton's method including line search.
		The required iteration numbers are displayed in Table~\ref{tab: Newtondata}.
		Only in the first refinement step, a line search is required. In all other refinement steps a few Newton iterations are necessary until the method is converging.
		There are about three for the Density model and mostly five or six for the Boussinesq model.
	}

        \newpage
	\section{Conclusions}
	\label{Sec; Concl}
	In this work, we first proposed a model for coupling the Navier-Stokes 
	equations with a temperature equation in which the density and viscosity 
	both are temperature-dependent. This results into a highly nonlinear problem. 
	We formulated a monolithic setting and applied multigoal-oriented 
	a posteriori error estimation up to six goal functionals (quantities of interest). 
	Therein, adjoint sensitivities measures for the error estimator 
	are obtained with the help of the dual-weighted residual method.
	Then, motivated from laser material processing, we conducted three numerical tests.
	Therein, due to the highly nonlinear behavior, our Newton solver would 
	only work with a priori local mesh refinement in certain regions of interest. 
	Moreover, we investigated error reductions and effectivity indices. According 
	to the highly nonlinear behavior and various parameters variations, 
	we obtained satisfying results
	for both error reductions and effectivity indices.
	
	\section*{Acknowledgments}
	This work has been supported by the Cluster of Excellence PhoenixD (EXC 2122, Project ID 390833453).
	Furthermore, Bernhard Endtmayer is funded by an Humboldt Postdoctoral Fellowship and greatly acknowledges this. Thomas Wick thanks \v{S}\'arka Ne\v{c}asov\'a (Czech Academy of Sciences) for discussions in the course of this work on the Navier-Stokes equations with material dependent coefficients.
	
	
	\bibliographystyle{abbrv}
	\bibliography{lit}

\end{document}